\documentclass[final]{siamart190516}
\usepackage{t1enc}
\usepackage{amssymb}
\usepackage{amsmath}
\usepackage{latexsym}
\usepackage{stmaryrd}
\usepackage{wasysym}
\usepackage{booktabs}
\usepackage{multicol}
\usepackage{multirow}
\usepackage{graphicx}
\usepackage{xcolor}
\usepackage{paralist}
\usepackage{arydshln}
\usepackage[mathcal]{euscript}
\usepackage{algorithmic}
\usepackage{pbox}
\usepackage{nicefrac}

\newcommand{\assgn}{\ensuremath\mathrel{\mathop:}=}
\newcommand{\sgn}{\operatorname{sgn}}

\newcommand{\diag}{\operatorname{diag}}
\newcommand{\Real}{\operatorname{Re}}
\newcommand{\Imag}{\operatorname{Im}}
\newcommand{\rv}{\textbf{r}}
\newcommand{\nbr}{\operatorname{\mathit{n}}(\hbox{\textbf{r}})}
\newcommand{\kv}{\hbox{\textbf{k}}}
\newcommand\hsumpart[2]{#1_{L'\!, \mkern1mu a,t'}^{*}\thinspace
  T_{L'\!,\mkern1mu L;a}^{[#1#2]} \thinspace
  #2_{L,a,t}^{}}
\newcommand{\odg}[1]{\mathcal{O}(#1)}
\title{The LAPW method with eigendecomposition based on the
  Hari--Zimmermann generalized hyperbolic SVD%
\thanks{This work has been supported in part by Croatian Science
  Foundation under the project IP--2014--09--3670, and also in part by
  a bilateral research project ``Optimization of material science
  algorithms on hybrid HPC platforms'' funded by the Croatian Ministry
  of Science and Education (MZO) and the German Academic Exchange
  Service (DAAD).}}%
\author{Sanja Singer\thanks{University of Zagreb, Faculty of
    Mechanical Engineering and Naval Architecture, I.~Lu\v{c}i\'{c}a
    5, 10000 Zagreb, Croatia, (ssinger@fsb.hr).}
  \and Edoardo Di Napoli\thanks{Forschungszentrum J\"{u}lich,
    J\"{u}lich Supercomputing Centre, Wilhelm--Johnen-Stra\ss e,
    J\"{u}lich, 52425, Germany and RWTH Aachen University, AICES,
    Schinkelstra\ss e 2, Aachen, 52062, Germany,
    (e.di.napoli@fz-juelich.de, dinapoli@aices.rwth-aachen.de).}
  \and Vedran Novakovi\'{c}\thanks{Completed a majority of his part of
    the research while being affiliated to Universidad Jaime I,
    Av.~Vicent Sos Baynat, 12071 Castell\'{o}n de la Plana, Spain,
    (novakoni@uji.es).}
  \and Gayatri \v{C}aklovi\'{c}\thanks{Ph.D.~student,
    Forschungszentrum J\"{u}lich, J\"{u}lich Supercomputing Centre,
    Wilhelm--Johnen-Stra\ss e, J\"{u}lich, 52425, Germany,
    (g.caklovic@fz-juelich.de).}}
\begin{document}
\maketitle
\renewcommand{\thefootnote}{\fnsymbol{footnote}}
\begin{abstract}
  In this paper we propose an accurate, highly parallel algorithm for
  the generalized eigendecomposition of a matrix pair $(H, S)$, given
  in a factored form $(F^{\ast} J F, G^{\ast} G)$.  Matrices $H$ and
  $S$ are generally complex and Hermitian, and $S$ is positive
  definite.  This type of matrices emerges from the representation of
  the Hamiltonian of a quantum mechanical system in terms of an
  overcomplete set of basis functions.  This expansion is part of a
  class of models within the broad field of Density Functional Theory,
  which is considered the golden standard in condensed matter physics.
  The overall algorithm consists of four phases, the second and the
  fourth being optional, where the two last phases are computation of
  the generalized hyperbolic SVD of a complex matrix pair $(F,G)$,
  according to a given matrix $J$ defining the hyperbolic scalar
  product.  If $J = I$, then these two phases compute the GSVD in
  parallel very accurately and efficiently.
\end{abstract}
\begin{keywords}
  LAPW method,
  generalized eigendecomposition,
  generalized (hyperbolic) singular value decomposition,
  hyperbolic QR factorization
\end{keywords}
\begin{AMS}
  65F15, 65F25, 65Y05, 65Z05
\end{AMS}
\pagestyle{myheadings}
\thispagestyle{plain}
%
%
\section{Introduction}\label{s:1}
%
%
Density Functional Theory (DFT) is the Standard Model at the base of
simulations in condensed matter physics.  At the center of most DFT
simulations lays the initialization of the Hamiltonian matrix $H$ and
its diagonalization.  In many DFT methods the form and size of the
Hamiltonian depends on the choice of the set of basis functions used
to expand the atomic orbitals.  When such a basis set is not
orthonormal, a Hermitian positive definite overlap matrix $S$ has to
be computed and diagonalized simultaneously with $H$; this pair of
matrices $(H, S)$ define a generalized Hermitian eigenproblem (or
eigenpencil, in short).  In a subset of all DFT methods labeled as
LAPW, the entries of both $H$ and $S$ are represented as multiple sums
and products of smaller matrices with specific properties.  We show
how to exploit this peculiar representation to solve the generalized
eigenvalue problem without explicitly assembling the $H$ and $S$
matrices.  Our alternative method solves the eigenpencil using a
cascade of phases ending with the Hari--Zimmermann algorithm for a
generalized hyperbolic SVD.  We demonstrate the scalability of a
shared memory version of this method on a number of test cases
extracted from concrete DFT simulations.  If the matrix $S$ is
ill-conditioned, our method has the additional benefit of providing
enhanced accuracy while avoiding the failure-prone Cholesky
factorization.

The birth of DFT is marked by two fundamental articles by the Nobel
prize winner Walter Kohn and his collaborator Lu J.~Sham and Pierre
Hohenberg~\cite{Hohenberg-64,Kohn-65}.  DFT provides an approach to
the theory of electronic structure that is alternative to the solution
of the Schr\"odinger equation.  While in the latter the emphasis is on
a many-electron wave function describing the dynamics of electrons in
a multi-atomic system, in DFT the electron density distribution $\nbr$
plays a central role.  Besides providing a complementary perspective,
DFT has made possible the simulation of much larger systems than the
conventional multi-particle wave function methods.  Depending on the
specific DFT method, computing complexity scales at most with the cube
of the number of atoms, with ongoing progress towards bringing it down
to linear scaling.

Despite being a general theory, DFT can be realized in as many flavors
as are the sets of basis functions one can choose from.  Two widely
spread classes of basis functions build on the simplicity of plane
waves to build more complex and rich sets of basis functions, namely
Projected Plane Waves (PAW)~\cite{Rostgaard-09} and Linearized
Augmented Plane Waves (LAPW)~\cite{Singh-06}.  The complexity of these
sets lays in that they are made up of non-orthogonal basis functions.
In the case of LAPW, the set of functions is also overcomplete.  The
consequence of non-orthogonality is that the matrix $S$, whose entries
are the scalar products among all the basis functions of any given
finite size set, is usually dense.  In the particular case of LAPW
methods such matrix is positive definite but could have few singular
values quite close to zero.  This potential problem is due to the
overcompleteness of the basis set and tends to worsen as the number of 
atoms increases since the number of basis functions grows linearly
with it.

In DFT methods the dynamics of the quantum systems is described by a
Hamiltonian operator.  In practice, the Hamiltonian is translated into
a Hermitian matrix $H$ whose size and structure depends on the
specific DFT method.  This is because the matrix $H$ is the result of
the projection of the Hamiltonian operator over the finite set of
basis functions of the given method.  In the LAPW method, the
mathematical form of the functions leads to an expression for both $H$
and $S$ in terms of a sum of smaller matrices over all possible
atoms $N_A$,
\begin{equation}
  \begin{aligned}
    H & = \sum_{a = 1}^{N_A} (A_a^{\ast} T_a^{[A A]} A_a^{} + A_a^{\ast} T_a^{[A B]} B_a^{}
    + B_a^{\ast} T_a^{[B A]} A_a^{} + B_a^{\ast} T_a^{[B B]} B_a^{}),\\
    S & = \sum_{a = 1}^{N_A} (A_a^{\ast} A_a^{} + B_a^{\ast} U_a^{\ast} U_a^{} B_a^{}),
  \end{aligned}
  \label{1.1}
\end{equation}
where $A_a, B_a \in \mathbb{C}^{N_L \times N_G}$, with $N_G$ and $N_L$
($N_G \geq N_L$) being the size of the basis set and the total number
of angular momentum states, respectively.  The remaining matrices
in~(\ref{1.1}) are complex, square, of order $N_L$, with some
additional properties.  Matrices $U_a^{}$ are real and diagonal, for
all $T_a^{[A B]}$ holds $(T_a^{[A B]})^{\ast} = T_a^{[B A]}$, while
$T_a^{[A A]}$ and $T_a^{[B B]}$ are Hermitian.  Except for $U_a$, the
other matrices are in general dense, and can have a range of sizes
dictated by the constants $N_A$, $N_G$, and $N_L$ (see
section~\ref{s:2} for some their typical range).  Despite the
formulation above could lend itself to computation through specialized
middleware libraries such as the Basic Linear Algebra Subprograms
(BLAS), the standard approach followed by most code developers was one
based on minimizing memory footprint and FLOP
count~\cite{Canning-00,Kurz-00}.

Recently, an alternative method for the assembly of the matrices $H$
and $S$ was presented
in~\cite{DiNapoli-Peise-Hrywniak-Bientinesi-2017} and further
developed
in~\cite{Davidovic-FabregatTraver-Hoehnerbach-DiNapoli-2018}.  In
their work~\cite{DiNapoli-Peise-Hrywniak-Bientinesi-2017}, Di~Napoli
et al.\ consolidate the underlying matrix structure of the operations
and proceed to encapsulate them in terms of the level 3 BLAS kernels.
For instance, to maximize the arithmetic intensity of the computation,
matrix $H$ is written as $H = H_{AA} + H_{AB+BA+BB}$, with
\begin{displaymath}
  H_{AA} = \sum_{a = 1}^{N_A} A_a^{\ast} T_a^{[A A]} A_a^{}, \quad
  H_{AB+BA+BB} = \sum_{a = 1}^{N_A} \left(B_a^{\ast} Z_a^{} + Z_a^{\ast} B_a^{} \right),
\end{displaymath}
where $Z_a^{} = T_a^{[B A]} A_a^{} + \nicefrac{1}{2} T_a^{[B B]} B_a^{}$.

Each of the $Z_a$ and $B_a$ matrices is then packed in memory in two
consecutive 2-dimensional arrays $Z_*$ and $B_*$, respectively.  In
the end, the sum $H_{AB+BA+BB}$ is computed by just two
\texttt{ZHER2K} BLAS subroutines.  A similar procedure holds for the
matrix $S$.  Once assembled the algebraic dense generalized
eigenproblem is solved by standard methods.  A Cholesky factorization
$LL^{\ast} = S$ is used to reduce the problem to standard form
$A \leftarrow L^{-1}AL^{-\ast}$.  In turn, the standard problem is
solved by a dense direct algorithm such as MRRR~\cite{Dhillon-04}
provided by the LAPACK library~\cite{Anderson-et-al-99}, or an
iterative eigensolver specialized for DFT computation (e.g., the ChASE
library~\cite{Winkelmann-19}).  When the assembled $S$ matrix is
ill-conditioned, as it may happen for quantum systems with a large
number of atoms ($> 100$), the Cholesky factorization may fail and
makes it practically very hard to solve the corresponding generalized
eigenproblem.  This issue is typically solved by the practitioners by
modifying the mathematical model so to avoid increasing the basis set,
which is the source of an ill-condioned $S$, but thus compromising on
the robustness of the DFT approach.

In this work, we propose an numerical method alternative to the
physics-based approach for solving the eigenpencil $(H, S)$ without
forming the matrices explicitly.  The core of the method is based on
the generalized hyperbolic singular value decomposition
(GHSVD)~\cite{Bojanczyk-2003}.  Not only such a method solves for the
eigenproblem directly without assembling $H$ and $S$, but also could
give more accurate results when $S$ is nearly singular.  This is
possible since the GHSVD decomposition acts directly on the
multiplying factors making up $S$, conceivably reducing the
singularity down to the square root of the condition number of $S$.
As a surplus, if $J$, the matrix of the hyperbolic scalar product, is
equal to the identity, the GHSVD reduces to the generalized SVD
(GSVD), which is computed very efficiently in parallel.

The paper is subdivided into eight sections.  In section~\ref{s:2}, we
present in more detail the physics of the problem and the mathematical
model leading to the expression~(\ref{1.1}) for $H$ and $S$.
Section~\ref{s:3} is devoted to formulating the problem in precise
algebraic terms, and provides an overview of our algorithm.  The next
four sections deal with the four phases of the algorithm, where the
first three of them belong to the algorithm proper, and the fourth one
completes the computation of the GHSVD and is unrelated to the
underlying mathematical physics origin of the problem.  Since each
phase is an algorithm and a reusable software contribution in its own
right, at the end of each section we present the numerical results and
the parallelization techniques applied.  The paper concludes with a
note on the related future work in section~\ref{s:8}.
%
%
\section{The $H$ and $S$ matrices in LAPW methods}\label{s:2}
%
%
At the core of DFT are a set of equations, called Kohn--Sham
equations, that have to be solved for each of the single particle wave
function $\psi_i$
\begin{equation}
  \hat{H}_{\rm KS}\, \psi_i(\rv)
  = \left[ -\frac{\hbar^2}{2m_e}\nabla_{\rv}^2
    + V[\nbr](\rv) \right] \psi_i(\rv)
  = \epsilon_i \psi_i(\rv),
  \quad i = 1, \ldots, N_e.
  \label{2.1}
\end{equation}
The peculiarity of these equations is that the Hamiltonian operator
$\hat{H}_{\rm KS}$ depends implicitly on all the $\psi$ through the
charge density function $\nbr$, which makes the entire set of
Kohn--Sham equation strongly coupled and non-linear.  In particular
the function $\nbr$ is the sum of the squares of all $\psi$ up to the
total number of electrons $N_e$ in any given quantum system
\begin{equation}
  \nbr = \sum_{i=1}^{N_e} |\psi_i(\rv)|^2.
  \label{2.2}
\end{equation}
Because the equations~(\ref{2.1}) are non-linearly coupled, they can
be solved only self-consistently: one starts from a reasonable guess
for the charge density $\nbr_{\rm start}$, computes the potential
$V[\nbr]$, and solves~(\ref{2.1}).  The resulting functions $\psi_i$
and values $\epsilon_i$ are then used to compute a new density as
in~(\ref{2.2}), which is compared to the starting one.  If the two
densities do not match, the self-consistent loop is repeated with a
new mixed charge density.  The loop stops only when the new and the
old density agree up to some defined constant.

So far we have described the general setup.  There are many methods
that translate this setup into an algorithm, and this is where the
various ``flavors'' of DFT differ.  The first difference is in the
choice of the set of functions $\varphi_t$ used to expand every one
particle wave function $\psi_i$
\begin{equation}
  \psi_i(\rv) = \sum_{t=1}^{N_G} c_{t,i} \, \varphi_t(\rv).
  \label{2.3}
\end{equation}
In the LAPW method~\cite{Jansen-84,Wimmer-81}, the configuration space
where the atomic cells are defined is divided in two disjoint areas
where the wave functions have distinct symmetries: close to the atomic
nuclei, solutions tend to be spherically symmetric and strongly
varying, while further away from the nuclei, they can be approximated
as uniformly oscillating.  The qualitative structure of the solution
leads to a space composed of non-overlapping spheres---called muffin
tins (MT)---separated by interstitial (INT) areas.  The complete set
of basis functions $\varphi_t$ are given by a piece-wise definition
for each of the $N_A$ atoms and relative surrounding regions.
\begin{equation}
  \varphi_t(\rv)
  = \begin{cases}
    \displaystyle
    \sum_{l=0}^{l_{\max}}
    \sum_{m=-l}^{l}
    \left [ A_{(l,m),a,t} u_{l,a}(r)
      + B_{(l,m),a,t} \dot{u}_{l,a}(r) \right]
    Y_{l,m} (\hat{\bf r}_a), & a^{\rm th}\ \text{MT}\\[12pt]    
    \displaystyle \frac{1}{\sqrt{\Omega}}
    \exp(\mathrm{i} {\bf k}_t \cdot \rv), & \text{INT}.
  \end{cases}
  \label{2.4}
\end{equation}
In the MT spheres, each basis function depends on specialized radial
functions $u_{l,a}$, their derivatives $\dot u_{l,a}$ and the spherical
harmonics $Y_{l,m}$; the former only depend on the distance {\bf r}
from the MT center, while the latter form a complete basis on the unit
sphere defined by $\hat{\bf r} = \nicefrac{\rv}{|\rv|}$ and so depends
solely on the MT spherical angles.  Despite being piece-wise
functions, $\varphi_t$ must be continuous and differentiable for each
index $t$ and each atomic index $a$.  The coefficients
$A_{l,m,a}, B_{l,m,a} \in \mathbb{C}$ are set to guarantee that
$\varphi_t \in C^1$ for each of the values of the indices
$L \equiv (l, m)$ and $a$.  The variable $t$ ranges over the size of
the plane wave functions set in INT, and is used to label the vector
$\kv_t$ living in the space reciprocal to ${\rv}$.  As such, the
momentum $\kv_t$ characterizes the specific wave function entering in
the basis set.  The total size of the basis set is determined by
setting a cutoff value ${\bf K}_{\max} \geq \kv_t$.

When one substitutes the expansion of $\psi_i$~(\ref{2.3})
in~(\ref{2.1}), the Kohn--Sham equations become an algebraic
generalized eigenvalue problem that needs to be solved for the
$N_G$-tuples of coefficients
$c_i = (c_{1,i}, \ldots, c_{N_G,i})^T$
\begin{displaymath}
  \sum_{t=1}^{N_G} (H)_{t'\!,\mkern1mu t} \, c_{t,i}
  = \epsilon_i \sum_{t=1}^{N_G} (S)_{t'\!,\mkern1mu t} \, c_{t,i}.
\end{displaymath}
The complexity of the LAPW basis set is transferred to the definition
of the entries of the Hamiltonian and overlap matrices, respectively
$H$ and $S$, given by
\begin{displaymath}
  (H)_{t'\!,\mkern1mu t} = \sum_a \iint \varphi^{\ast}_{t'}(\rv) \hat{H}_{\rm KS}^{} \,
  \varphi_t^{}(\rv) \, {\rm d}\rv, \quad
  (S)_{t'\!,\mkern1mu t} =  \sum_a \iint \varphi^{\ast}_{t'}(\rv) 
  \varphi_t^{}(\rv) \, {\rm d}\rv.
\end{displaymath}
By substituting explicitly the functions $\varphi_t$ of
equation~(\ref{2.4}) and computing the integrals, one ends up
with the following expressions for $H$ and $S$:
\begin{align}
  (H)_{t'\!,\mkern1mu t} & = \sum_a \sum_{L'\!,\mkern1mu L}
  \left(\hsumpart AA\right)
  + \left(\hsumpart AB\right) \label{2.5} \\
  & \qquad\qquad {} + \left(\hsumpart BA\right) + \left(\hsumpart BB\right),
  \nonumber \\
  (S)_{t'\!,\mkern1mu t}
  & = \sum_a \!\! \sum_{L=(l,m)} A_{L,a,t'}^{*} A_{L,a,t}^{}
  + B_{L,a,t'}^{*} B_{L,a,t}^{} \| \dot u_{l,a}^{} \|^2.
  \label{2.6}
\end{align}
The new matrices
$T_{L'\!,\mkern1mu L;a}^{[\dots]}\in \mathbb C^{N_L\times N_L}$ are
dense and their computation involves multiple integrals between the
radial basis functions $u_{l,a}$ and the non-spherical part of the
potential $V$ multiplied by Gaunt coefficients (for details
see~\cite{Kurz-00}\cite[Ch.~5]{Singh-06}\cite[App.]{DiNapoli-Peise-Hrywniak-Bientinesi-2017}).
As can be seen by simple inspection,
equations~(\ref{2.5})--(\ref{2.6}) are equivalent to
equations~(\ref{1.1}): while the former are written with all indices
explicit, the latter have a subset of them implicit which highlights
their matrix form.

\looseness=-1
We conclude with a small excursus on the structure of the
self-consistent loop and its computational cost.  In the first step, a
starting charge density $\nbr_{\rm start}$ is used to compute the
Kohn--Sham Hamiltonian $H_{\rm KS}$.  In a second step the set of
basis functions is set up and the set of $A,B$ coefficients is
derived.  Then, the Hamiltonian $H$ and overlap $S$ matrices are
initialized, followed by the fourth step when the generalized
eigenvalue problems $H c = \epsilon\ S c$ is solved numerically to
return the eigenpairs $(C, {\rm diag}(\epsilon))$.  Finally a new
charge density $\nbr$ is computed and convergence is checked before
starting a new loop.  Out of all the steps above, initializing $H$ and
$S$ and solving the eigenproblem accounts for more than 80\%\ of CPU
time.  Having cubic complexity $\odg{N_G^3}$, the eigenproblem
solution is usually considered the most expensive of the two.
It turns out that generating the matrices may be as expensive.  If
$N_A$ and $N_L$, respectively, are the range of the summations
$\sum_a$ and $\sum_L$, then it can be shown that equations~(\ref{2.6})
and (\ref{2.5}) have complexity $\odg{N_A \cdot N_L \cdot N_G^2}$ and
$\odg{N_A \cdot N_L \cdot N_G \cdot(N_L + N_G)}$.  A typical
simulation uses approximately $N_G$ basis functions, with $N_G$
ranging from about $50 \cdot N_A$ to about $80 \cdot N_A$, and an
angular momentum $l_\mathrm{max} \leq 10$, which results in
$N_L = (l_{\max} + 1)^2 \leq 121$.  It follows that the factor
$N_A \cdot N_L$ is roughly of the same order of magnitude as $N_G$ so
that the generation of $H$ and $S$ also displays cubic complexity
$\odg{N_G^3}$.  In practice, the constants above have values in the
following orders of magnitude:
$N_A = \mathcal{O}(100)$,
$N_G = \mathcal{O}(1000)$--$\mathcal{O}(10000)$, and
$N_L = \mathcal{O}(100)$.
%
%
\section{Problem formulation}\label{s:3}
%
%
Our intention is to keep the matrices $H$ and $S$ in their factored
form given in (\ref{1.1}).  The core of the process, Phase~3, is a
one-sided Jacobi-like method for the implicit diagonalization, that
computes a hyperbolic analog of the generalized SVD.
\begin{definition}
  For the given matrices $F \in \mathbb{C}^{m \times n}$, $m \geq n$,
  $J \in \mathbb{R}^{m \times m}$, $J = \diag(\pm 1)$, and
  $G \in \mathbb{C}^{p \times n}$, where $G$ is of full column rank,
  there exist a $J$-unitary matrix $U \in \mathbb{C}^{m \times m}$
  (i.e., $U^{\ast} J U = J$), a unitary matrix
  $V \in \mathbb{C}^{p \times p}$, and a nonsingular matrix
  $X \in \mathbb{C}^{n \times n}$, such that
  \begin{equation}
    F = U \Sigma_F X, \quad
    G = V \Sigma_G X, \quad
    \Sigma_F \in \mathbb{R}^{m \times n}, \quad
    \Sigma_G \in \mathbb{R}^{p \times n}.
    \label{3.1}
  \end{equation}
  The elements of $\Sigma_F$ and $\Sigma_G$ are zeros, except for the
  diagonal entries, which are real and nonnegative.  Furthermore,
  $\Sigma_F$ and $\Sigma_G$ satisfy
  $\Sigma_F^T \Sigma_F^{} + \Sigma_G^T \Sigma_G^{} = I$.  The ratios
  $\Sigma_{ii} \assgn (\Sigma_F)_{ii}/(\Sigma_G)_{ii}$ are
  called the generalized hyperbolic singular values of the pair
  $(F, G)$.  If the pair $(F, G)$ is real, then all matrices in
  (\ref{3.1}) are real.
\end{definition}
We choose to define the generalized hyperbolic SVD (GHSVD) only if the
matrix $G$ is of full column rank.  This implies $p \geq n$, and there
is no need to mention this in the definition.  In the case of full
column rank $G$, the matrix $S \assgn G^{\ast}G$ is positive definite,
and the matrix pair $(H, S)$, where  $H \assgn F^{\ast} J F$, is
Hermitian and definite, so it can be simultaneously diagonalized by
congruences (see, for example,~\cite{Parlett-98a}).

If the GHSVD is computed as in (\ref{3.1}), then the generalized
eigenvalues and eigenvectors of $(H, S)$ are easily retrieved, since
\begin{align*}
  H & = F_{}^{\ast} J F
    = X_{}^{\ast} \Sigma_F^{\ast} U_{}^{\ast} J U \Sigma_F^{} X
    = X_{}^{\ast} \Sigma_F^{\ast} J \Sigma_F^{} X
    \assgn X_{}^{\ast} \Lambda_F^{} X, \\
  S & = G_{}^{\ast} G
    = X_{}^{\ast} \Sigma_G^{\ast} V_{}^{\ast} V \Sigma_G^{} X
    = X_{}^{\ast} \Sigma_G^{\ast} \Sigma_G^{} X
    \assgn X_{}^{\ast} \Lambda_G^{} X.
\end{align*}
Substituting $X^{\ast} = S^{} X^{-1} \Lambda_G^{-1}$ in the expression
for $H$ above, we get
\begin{equation}
  H^{} Z^{} = S^{} Z^{} \Lambda;\quad Z \assgn X^{-1}, \quad
  \Lambda \assgn \Lambda_G^{-1} \Lambda_F^{}.
  \label{3.2}
\end{equation}
Thus, from (\ref{3.2}), the generalized eigenvalues $\diag(\Lambda)$
of the matrix pair are the squared generalized hyperbolic singular
values, with the signs taken from the corresponding diagonal elements
in $J$, i.e., $\diag(\Sigma^T J \Sigma)$, and the matrix of the
generalized eigenvectors $Z$ is the inverse of the matrix $X$ of the
right generalized singular vectors.  For theoretical purposes it can
be assumed that $\diag(\Lambda)$ is sorted descendingly, though for
simplicity it is not the case in our implementation.

An approach that uses the SVD on a matrix factor, instead of the
eigendecomposition on the multiplied factors, usually computes small
eigenvalues more accurately.

In the first phase of the algorithm, we transform the initial problem
by assembling the Hermitian matrices $T_a$,
\begin{equation}
  T_a =
  \begin{bmatrix}
    T_a^{[A A]} & T_a^{[A B]}\\
    T_a^{[B A]} & T_a^{[B B]}
  \end{bmatrix}
  \label{3.3}
\end{equation}
and factoring them into a form suitable for the GHSVD computation.

After that, we are left with two tall matrices, which have, in our
test examples, between $2$ and $22$ times more rows than columns.
Since the Jacobi-like SVD algorithms are more efficient if the factors
are square, in the second (optional) phase we could preprocess the
factors: $F$ by the hyperbolic QR factorization
(see~\cite{SingerSanja-2006}), and $G$ by the tall-and-skinny QR
factorization (\texttt{ZGEQR} routine from LAPACK), into square ones.

The third phase is a complex version of the implicit Hari--Zimmermann
method---a modification of the one-sided real method presented
in~\cite{Novakovic-SingerSanja-SingerSasa-2015}.  The complex
transformations, for the two-sided method, were derived by Vjeran Hari
in his PhD thesis~\cite{Hari-84}.
\subsection{Overview of the algorithm}\label{ss:3.1}
The sequence of phases of our algorithm is:
\begin{compactenum}[1.]
\item The problem is expressed as
  $H = F_0^{\ast} \diag(T_1, \ldots, T_{N_A}) F_0^{}$,
  $S = \widetilde{G}^{\ast} \widetilde{G}$, the matrices $F_0$ and
  $\widetilde{G}$ are assembled, and the matrices $T_a$, formed from
  $T_a^{[A A]}$, $T_a^{[A B]}$, $T_a^{[B A]}$, and $T_a^{[B B]}$, are
  simultaneously factored by the Hermitian indefinite factorization
  with complete pivoting, reformulating $H$ as
  $H = \widetilde{F}^{\ast} \widetilde{J} \widetilde{F}^{}$, with
  $\widetilde{J} = \diag(\pm 1)$.
\item Optionally, the tall-and-skinny matrices $\widetilde{F}$ and
  $\widetilde{G}$ are shortened: $\widetilde{F}$ by the indefinite,
  $\widetilde{J}$-QR factorization to obtain the square factor $F$ and
  a new signature matrix $J$, and $\widetilde{G}$ by the QR
  factorization to obtain the square factor $G$.
\item The GEVD of $(H, S)$ is computed by the $J$-GHSVD of $(F, G)$
  (or the $\widetilde{J}$-GHSVD of $(\widetilde{F}, \widetilde{G})$ if
  the Phase~2 is skipped) by the implicit Hari--Zimmermann method.
\item Optionally, the GHSVD process is formally completed by
  explicitly computing  the right generalized singular vectors $X$
  from the generalized eigenvector matrix $Z$.
\end{compactenum}
\subsection{Testing environment and data}\label{ss:3.2}
\looseness=-1
The testing environment consists of a node with an Intel Xeon Phi 7210
CPU, running at 1.3~GHz with Turbo Boost turned off, in Quadrant
cluster mode with 96~GiB of RAM and 16~GiB of flat-mode MCDRAM, under
64-bit CentOS Linux 7 with the Intel compilers (Fortran, C) and Math
Kernel Library (MKL) version 19.0.5.281, and GNU Fortran 8.3.1 for the
error testing.

The software code, freely available in
\url{https://github.com/venovako/FLAPWxHZ} repository, of all the
phases presented in this paper is written mostly in Fortran, with some
auxiliary parts in C, while the parallelization relies on the OpenMP
constructs.

\looseness=-1
The phases are meant to be run in a sequence, where each phase is
executed as a separate process with several OpenMP threads.  Since the
modern compute nodes generally have enough memory to hold all required
data, the algorithms are implemented for the shared memory, but at
least the algorithms for the Phases~1, 3, and 4 can be transformed
into distributed-memory ones, should the volume of data so require.

In testing it was established that each thread should be bound to its
own physical CPU core, with \texttt{OMP\_PROC\_BIND=SPREAD} placement
policy.  The double precision and the double-complex BLAS and LAPACK
routines were provided by the thread-parallel MKL, but with only one
(i.e., the calling) thread allowed per call, except for the
\texttt{ZSWAP}, \texttt{ZROT}, and \texttt{ZGEQR} routines in Phase~2,
\texttt{ZGETC2} routine in Phase~4, and \texttt{ZGEMM},
\texttt{ZHERK}, \texttt{ZHEGV}, and \texttt{ZHEGVD} routines in
subsection~\ref{sss:6.5.4}, where the MKL was allowed to use the
test's upper limit on the number of threads.  The nested parallelism
is therefore possible but not required in our code.  Hyper~Threading
was enabled but not explicitly utilized, though nothing precludes a
possibility that on a different architecture the BLAS or LAPACK calls
could benefit from some form of intra-core symmetric multi-threading.
A refined thread placement policy
\texttt{OMP\_PROC\_BIND=SPREAD,CLOSE} might then allow better reuse of
data in the cache levels shared among the threads of a core.

Apart from the maximal number of threads set to the number of CPU
cores in a node, the tests were also performed with half that number,
to assess the effects on the computational time of the larger block
sizes and the availability of the whole L2 data cache (1~MiB, shared
among two cores) to a thread.  The algorithms do not constrain the
number of threads in principle, but are not intended to be used
single-threadedly.

Our main test node has 64 cores, but a subset of the tests were
repeated on a faster JUWELS~\cite{JUWELS} node, with two Intel Xeon
Platinum 8168 CPUs, running at 2.7~GHz with 1~MiB L2 cache per each of
$2\times 24=48$ cores, with a similar software setup, for a comparison
of the GHSVD and the generalized eigendecomposition approaches (see
subsection~\ref{sss:6.5.4}).  When the results obtained on JUWELS are
shown, the test's number of threads is emphasized (e.g., \emph{48\/}),
to distinguish them from the main results.

Another hardware feature targeted is the SIMD vectorization: each core
of both machines has a private L1 data cache of 32~kiB with a line
size of 64~B, and equally wide (e.g., 8 double precision
floating-point numbers) vector registers upon which a subset of
AVX-512 instructions is capable to operate in the SIMD fashion.  The
vectorization is employed both implicitly, by aligning the data to the
cache line size whenever possible and instructing the compiler to
vectorize the loops, and semi-explicitly, as will be described in the
following sections.  The code is parametrized by the maximal SIMD
length (i.e., the number of 8~B lanes in the widest vector register
type) $\mathsf{v}$, and it vectorizes successfully on other
architectures (e.g., on AVX2, with $\mathsf{v}=4$).

Under an assumption that the compiler-generated floating-point
reductions (e.g., those of the \texttt{SUM} Fortran intrinsic) obey
the same order of operations in each run, and due to the alignment
enforced as above, the algorithms should be considered conditionally
reproducible, in a sense that the multiple runs of the same
executables on the same data in the same environment should produce
bitwise-identical results.
\subsubsection{Datasets}\label{sss:3.2.1}
Each dataset under test contained all matrix inputs ($A_a^{}$,
$B_a^{}$, $U_a^{}$, $T_a^{[A A]}$, $T_a^{[B A]}$, $T_a^{[B B]}$) for a
single problem instance.  With 8 datasets from Table~\ref{tbl:ds} we
believe to have a representative coverage of the small-to-medium size
problems from practice.
\begin{table}[hbt]
  {\footnotesize\caption{The datasets under test. For A datasets
    $N_L = 121$, $N_A = 108$, and $m = 2 N_L N_A = 26136$, while for
    B datasets $N_L = 49$, $N_A = 512$, and $m = 2 N_L N_A = 50176$.
    Also, $n = N_G$.}
    \label{tbl:ds}
    \begin{center}
      \begin{tabular}{@{}ccc@{}}
        \toprule
        ID & \texttt{AuAg} & $n$ \\
        \toprule
        A1 & \texttt{2.5} & \hphantom{0}3275 \\
        \midrule
        A2 & \texttt{3.0} & \hphantom{0}5638 \\
        \bottomrule
      \end{tabular}
      \quad
      \begin{tabular}{@{}ccc@{}}
        \toprule
        ID & \texttt{AuAg} & $n$ \\
        \toprule
        A3 & \texttt{3.5} & \hphantom{0}8970 \\
        \midrule
        A4 & \texttt{4.0} & 13379 \\
        \bottomrule
      \end{tabular}
      \qquad
      \begin{tabular}{@{}ccc@{}}
        \toprule
        ID & \texttt{NaCl} & $n$ \\
        \toprule
        B1 & \texttt{2.5} & 2256 \\
        \midrule
        B2 & \texttt{3.0} & 3893 \\
        \bottomrule
      \end{tabular}
      \quad
      \begin{tabular}{@{}ccc@{}}
        \toprule
        ID & \texttt{NaCl} & $n$ \\
        \toprule
        B3 & \texttt{3.5} & 6217 \\
        \midrule
        B4 & \texttt{4.0} & 9273 \\
        \bottomrule
      \end{tabular}
  \end{center}}
\end{table}
As already mentioned in section~\ref{s:2}, the maximum value of the
momentum ${\bf K}_{\max}$ which appears as an index to the dataset
label (e.g., \texttt{AuAg\_2.5}) determines the size of the basis
functions set $N_G$.  This is why datasets with same label (e.g.,
\texttt{AuAg}) but different index (e.g., \texttt{2.5}
vs.\ \texttt{3.0}) have differing values for $N_G$.  In the following,
the datasets are referred to by their IDs.
%
%
\section{Phase~1 -- simultaneous factorizations of $T_a$ matrices}\label{s:4}
%
%
The goal of this section is to rewrite the problem (\ref{1.1}) in a
form suitable for GHSVD computation.
\subsection{Problem reformulation}\label{ss:4.1}
The first step is to write (\ref{1.1}) as
\begin{equation}
  H = \sum_{a=1}^{N_A} H_a^{\ast} T_a^{} H_a^{}, \quad
  S = \sum_{a=1}^{N_A} S_a^{\ast} S_a^{}, \quad
  H_a = \begin{bmatrix}
    A_a^{} \\
    B_a^{}
  \end{bmatrix}, \quad
  S_a = \begin{bmatrix}
    A_a^{} \\
    U_a^{} B_a^{}
  \end{bmatrix}.
  \label{4.1}
\end{equation}
Furthermore, the matrices in (\ref{4.1}) can be expressed as
\begin{equation}
  \begin{aligned}
    H & =
    \begin{bmatrix}
      H_1^{\ast} & \cdots & H_{N_A}^{\ast}
    \end{bmatrix}
    \diag(T_1^{}, \ldots, T_{N_A}^{})
    \begin{bmatrix}
      H_1^{\ast} & \cdots & H_{N_A}^{\ast}
    \end{bmatrix}^{\ast}
    \assgn F_0^{\ast} T F_0^{},\\
    S & =
    \begin{bmatrix}
      S_1^{\ast} & \cdots & S_{N_A}^{\ast}
    \end{bmatrix}
    \begin{bmatrix}
      S_1^{\ast} & \cdots & S_{N_A}^{\ast}
    \end{bmatrix}^{\ast}
    \assgn \widetilde{G}^{\ast} \widetilde{G}.
  \end{aligned}
  \label{4.2}
\end{equation}
In (\ref{4.2}), $\diag(T_1^{}, \ldots, T_{N_A}^{})$ stands for a
block-diagonal matrix with the prescribed diagonal blocks $T_a$,
$a = 1, \ldots, N_A$ from~(\ref{3.3}).  Newly defined matrices have
the following dimensions:
$H_a, S_a \in \mathbb{C}^{(2N_L) \times N_G}$,
$T_a \in \mathbb{C}^{(2N_L) \times (2N_L)}$,
$F_0, \widetilde{G} \in \mathbb{C}^{(2N_A N_L) \times N_G}$, and
$T \in \mathbb{C}^{(2N_A N_L) \times (2N_A N_L)}$.  From now on, let
$m \assgn 2N_A N_L$, and $n \assgn N_G$.

To efficiently exploit the structure of the problem, the matrix $T$
needs to be diagonal, with its diagonal elements equal to either $1$
or $-1$ (possibly with some zeros in the case of a singular $T$).
There is no theoretical obstacle to apply the simultaneous
($J$-)orthogonalization in the computation of the GHSVD on the
matrices $F_0$, $\widetilde{G}$, and $T$ implicitly, but the repeated
multiplication (in each reduction step) by $T$ is slow.  Therefore,
$T$ should be either factored concurrently, by using a modified
version of the Hermitian indefinite factorization of all $T_a$ blocks,
or diagonalized concurrently: all factorizations (or diagonalizations)
are independent of each other and can proceed in parallel.  Since the
diagonalization, compared to the Hermitian indefinite factorization,
is a slower process, our choice is to factor all the diagonal blocks
$T_a$.
\subsection{Hermitian indefinite factorization}\label{ss:4.2}
Each $T_a$ is factored by the algorithm described
in~\cite{Slapnicar-98}.  The algorithm for each $T_a$ consists of the
Hermitian indefinite factorization with a suitable
pivoting~\cite{Bunch-Parlett-71}, followed by the transformation of
the block-diagonal matrix.  Such factorization has the following form
\begin{equation}
  T_a = P_a^T M_a^{\ast} D_a^{} M_a^{} P_a^{},
  \label{4.3}
\end{equation}
where $P_a$ is a permutation (in the LAPACK sense), $M_a$ is upper
triangular, and $D_a$ is block-diagonal, with diagonal blocks of order
$1$ or $2$.

\looseness=-1
Then, $D_a$ is transformed into $\widehat{J}_a=\diag(\pm 1)$.  If
$D_a$ has a diagonal block of order $1$ at position $k$, then
$\widehat{J}_a$ stores the sign of this block in its $k$th diagonal
element, and the $k$th row of $M_a$ is scaled by
$|(D_a)_{kk}|^{1/2}$.  In the case of a (Hermitian) pivot block of
order $2$, this block is diagonalized by a Jacobi rotation
$R_k$, and the corresponding two rows of $M_a$ in (\ref{4.3}) are
multiplied by $R_k$.  Two transformations of the new diagonal elements
of $D_a$ are then performed, as above.  To speed-up the process, the
rotation and the scaling of two rows of $M_a$ are combined and then
applied as a single transformation.

The outer permutations $P_a$ are generated starting from the identity,
and stored as the partial permutations of the principal submatrices,
as in LAPACK, according to the pivoting of choice.  Since the matrices
$T_a$ are of a relatively small order, our choice is the complete
pivoting from~\cite{Bunch-Parlett-71}.
\subsection{Postprocessing}\label{ss:4.3}
After the factorization, with a postprocessing step we obtain
$T_a = \widehat{M}_a^{\ast} \widehat{J}_a^{} \widehat{M}_a^{}$, where
$\widehat{M}_a = M_a P_a$, and $\widehat{M}_a$ does not need to remain
triangular.

Finally, by applying an inner permutation $\widehat{P}_a$,
$\widehat{J}_a$ can be rearranged into a diagonal matrix
$\widetilde{J}_a$, where the positive signs precede the negative ones
on the diagonal.  This property of $\widetilde{J}_a$ matrices can be
exploited to speed-up computation of the hyperbolic scalar products
$x^{\ast} J^{} x^{}$ in the subsequent phases (see
subsections~\ref{sss:4.5.1}, \ref{ss:5.1}, and \ref{sss:6.4.2}).  Let
the whole factorization routine described thus far be called
\texttt{ZHEBPJ}.  Then,
$T_a = \widetilde{M}_a^{\ast} \widetilde{J}_a^{} \widetilde{M}_a^{}$,
$\widetilde{M}_a^{} = \widehat{P}_a \widehat{M}_a$, and $H_a$ is
multiplied by $\widetilde{M}_a$ from the left as
$\widetilde{H}_a = \widetilde{M}_a H_a$.

After such preprocessing, $H$ from (\ref{4.2}) is written as
\begin{displaymath}
  H = \widetilde{F}^{\ast} \widetilde{J} \widetilde{F}, \quad
  \widetilde{F}^{\ast}
  = [
    \widetilde{H}_1^{\ast}, \ldots, \widetilde{H}_{N_A}^{\ast}
  ], \quad
  \widetilde{J} = \diag(\widetilde{J}_1, \ldots, \widetilde{J}_{N_A}).
\end{displaymath}
In datasets A, each $\widetilde{J}_a$ has 3 positive and 239 negative
signs.  In datasets B, a non-con\-secutive half of $\widetilde{J_a}$
matrices are positive definite, and the others are negative definite.
\subsection{Implementation and testing}\label{ss:4.4}
The computational tasks for different indices $a$ are fully
independent, and are performed in parallel such that each thread is
responsible for one or more indices $a$, as indicated in the
pseudocode of Algorithm~\ref{alg:phase1}.
\begin{algorithm}
  \caption{A pseudocode for the Phase~1 algorithm.}
  \label{alg:phase1}
  \begin{algorithmic}
    \FORALL[an OpenMP parallel do]{atoms $a$, $1 \leq a \leq N_A$}
    \STATE{factorize $T_a = \widetilde{M}_a^{\ast} \widetilde{J}_a^{} \widetilde{M}_a^{}$;} \COMMENT{\texttt{ZHEBPJ} with BLAS level 1 and 2 routines}
    \STATE{multiply $\widetilde{H}_a = \widetilde{M}_a H_a$;} \COMMENT{$1$ \texttt{ZGEMM}, of a $2 N_L \times 2 N_L$ and a $2 N_L \times N_G$ matrix}
    \STATE{scale the rows of $B_a$ as $U_a B_a$;} \COMMENT{$N_L$ \texttt{ZDSCAL}s, each on a row of $N_G$ elements}
    \ENDFOR
  \end{algorithmic}
\end{algorithm}
Each thread, in turn, performs the three steps for its index $a$
sequentially, up to a possible usage of a parallel BLAS in the first
two steps.  The last step, i.e., computing $U_a B_a$ to assemble
$\widetilde{G}$, is a loop with the independent iterations, and could
be done in parallel, using the nested parallelism within each thread,
should $N_G$ be large enough and should also the newly spawned threads
for that loop have enough computational resources available to warrant
the overhead of the additional thread management.

Each thread is responsible for allocating (MCDRAM is not explicitly
used) and accessing the memory for the data it processes, so the data
locality is achievable whenever each NUMA node has enough storage.
The algorithm is thus viable in the heavily non-uniform memory
access settings, like the Intel Xeon Phi's SNC-4 mode.

In a distributed memory setting (e.g., using the MPI processes), the
assembling of $\widetilde{F}$, $\widetilde{J}$, and $\widetilde{G}$
can be done by assigning to each process a (not necessarily
contiguous) subrange of the iteration range of the for-all loop from
Algorithm~\ref{alg:phase1}, while inside the process all atoms
assigned to it are processed exactly as above, within an OpenMP
parallel-do loop.  The matrices $\widetilde{F}$, $\widetilde{J}$, and
$\widetilde{G}$ would then end up being distributed in the chunks
corresponding to the chosen subranges among the processes.
\subsubsection{Testing}\label{sss:4.4.1}
In Table~\ref{tbl:1} the average per-atom wall execution time of
Phase~1 is shown.
\begin{table}[hbt]
  {\footnotesize\caption{The average per-atom wall execution time
      (wtime) of Phase~1 with 32 and 64 threads.  Since the routine
      weights are rounded to the nearest per mil, their sum may not
      yield 100\%.  The first weight corresponds to \texttt{ZHEBPJ},
      the second one to \texttt{ZGEMM}, and the third one to
      \texttt{ZDSCAL}s step of Algorithm~\ref{alg:phase1}.}
    \label{tbl:1}
    \begin{center}
      \begin{tabular}{@{}ccccc@{}}
        \toprule
        \multirow[c]{2}{*}{ID}
        & \multicolumn{2}{c}{average wtime [s] per atom}
        & \multicolumn{2}{c}{routine weights \%:\%:\%} \\
        & 32 threads & 64 threads & 32 threads & 64 threads \\
        \toprule
        A1 & $0.243186$ & $0.277876$ & $71.9:26.2:\hphantom{0}1.9$ & $66.1:32.1:\hphantom{0}1.8$ \\
        \midrule
        A2 & $0.279318$ & $0.312057$ & $59.6:36.6:\hphantom{0}3.8$ & $53.8:41.4:\hphantom{0}4.9$ \\
        \midrule
        A3 & $0.345113$ & $0.409165$ & $47.9:46.4:\hphantom{0}5.7$ & $41.4:46.2:12.4$ \\
        \midrule
        A4 & $0.436803$ & $0.536776$ & $37.8:53.5:\hphantom{0}8.8$ & $31.7:50.9:17.5$ \\
        \midrule
        B1 & $0.023099$ & $0.027365$ & $56.3:38.4:\hphantom{0}5.3$ & $49.9:45.6:\hphantom{0}4.5$ \\
        \midrule
        B2 & $0.030430$ & $0.033248$ & $40.6:52.3:\hphantom{0}7.2$ & $38.6:54.8:\hphantom{0}6.7$ \\
        \midrule
        B3 & $0.045586$ & $0.060505$ & $25.6:65.1:\hphantom{0}9.3$ & $20.2:71.7:\hphantom{0}8.1$ \\
        \midrule
        B4 & $0.070669$ & $0.139536$ & $16.1:71.6:12.3$ & $\hphantom{0}8.3:80.9:10.8$ \\
        \bottomrule
      \end{tabular}
  \end{center}}
\end{table}
The results suggest that it is beneficial to have more L2 data cache
available per thread, as is the case with 32 threads overall.  In the
breakdown of the weights (i.e., percentages of time taken) of each
computational step it is confirmed that \texttt{ZGEMM} starts to
dominate the other computational steps of Algorithm~\ref{alg:phase1}
as the ratio $n/m$ increases.  It is a strong indication that even a
procedure more expensive than \texttt{ZHEBPJ}, such as a
diagonalization of $T_a$, may be applied on the datasets having a
square-like shape, without considerably degrading the relative
performance of Phase~1.

For a fully vectorized, cache-friendly alternative to applying
\texttt{ZDSCAL} with a non-unit stride in Algorithm~\ref{alg:phase1}
please refer to section~S.1 of the supplementary material.
\subsection{An alternative way forward}\label{ss:4.5}
After this phase has completed, one can proceed as described in the
rest of the paper, should the condition numbers of (the yet unformed)
matrices $H$ and $S$ be large enough to severely affect the accuracy
of a direct solution of the generalized Hermitian eigenproblem with
the pair $(H, S)$.

An alternative and more time-efficient way to proceed would be to
explicitly form $H$ and $S$.  For
$S^{} = \widetilde{G}^{\ast} \widetilde{G}^{}$, one \texttt{ZHERK}
call would suffice.  For
$H^{} = \widetilde{F}^{\ast} \widetilde{J}^{} \widetilde{F}^{}$, a
copy of $\widetilde{F}$ should be made, and that copy's rows should be
scaled in parallel by the diagonal elements of $\widetilde{J}$.  One
\texttt{ZGEMM} call on $\widetilde{F}^{\ast}$ and
$\widetilde{J}^{} \widetilde{F}^{}$ then completes the formation of
$H$.  After that, an efficient solver for the generalized Hermitian
eigenproblem can be employed on $(H, S)$, such as \texttt{ZHEGV} or
\texttt{ZHEGVD} from LAPACK, as shown in subsetion~\ref{sss:6.5.4}.
\subsubsection{Row scaling}\label{sss:4.5.1}
A cache-friendly implementation of the row scaling by $\widetilde{J}$
is to iterate sequentially over the rows of a fixed column $j$, and
change the sign of each element $\widetilde{F}_{ij}$ for which
$\widetilde{J}_{ii} = -1$, while the outer parallel-do loop iterates
over all column indices $j$.  However, that implementation can be
optimized further.

If $\widetilde{J}$ has its diagonal partitioned into (regularly or
irregularly sized) blocks of the same sign, then it can
be compactly encoded as a sequence of pairs $(i_{-},l)_k$, one for
each block of negative signs, where $i_{-}$ is the first index
belonging to a block $k$, and $l \ge 1$ is the block's length.  The
iteration over all rows and the conditional sign changes as above can
be replaced by iteration over all such blocks.  For each block,
iterate sequentially in the range of indices $i$ from $i_{-}$ to
$i_{-} + l - 1$, and change the signs unconditionally, thus
eliminating the conditional branching based on the sign of
$\widetilde{J}_{ii}$.

Such run-length-like encoding is employed in Phase~3, where it also
accelerates the hyperbolic dot products in the case where the positive
signs precede the negative ones on the diagonal of a sign matrix
(i.e., at most one negative block exists) given by \texttt{ZHEBPJ}
when forming the square factors for the inner Hari--Zimmermann method.
%
%
\section{Phase~2 -- optional $(J, I)$ URV factorization}\label{s:5}
%
%
The one-sided Jacobi-type algorithms are fastest if they work on
square matrices, since the column dot-products and updates are the
shortest possible.  If the square factors $F$, $G$, and the
corresponding $J$ of the matrix pair
$(\widetilde{F}^{\ast} \widetilde{J} \widetilde{F}, \widetilde{G}^{\ast} \widetilde{G})$,
can be found, instead of the rectangular factors $\widetilde{F}$,
$\widetilde{G}$, and the corresponding $\widetilde{J}$, we expect that
the overhead of such a shortening will be less than the computational
time saved by avoiding the rectangular factors.  To this end, matrix
$\widetilde{F}$ is shortened by using the hyperbolic QR factorization
(also called the JQR factorization), according to the given
$\widetilde{J}$:
\begin{equation}
  P_1^{} \widetilde{F}^{} P_2^{} = Q_F^{} F^{}; \quad
  \widetilde{Q}_F^{\ast} \widetilde{J}^{} \widetilde{Q}_F^{} = J^{}, \quad
  \widetilde{Q}_F^{} \assgn P_1^T Q_F^{},
  \label{5.1}
\end{equation}
where
$F \in \mathbb{C}^{n \times n}$ is block upper triangular with
diagonal blocks of order $1$ or $2$,
$J = \diag(\pm 1) \in \mathbb{R}^{n \times n}$ is the shortened
signature matrix, $P_2 \in \mathbb{R}^{n \times n}$ and
$P_1 \in \mathbb{R}^{m \times m}$ are the column and the row
permutation matrix, respectively.  Our application does not use
$Q_F^{} \in \mathbb{C}^{m \times n}$, so it is not explicitly formed.
From (\ref{5.1}) it holds
\begin{displaymath}
  P_2^T \widetilde{F}^{\ast} \widetilde{J}^{} \widetilde{F}^{} P_2^{}
  = F^{\ast} Q_F^{\ast} P_1^{} \widetilde{J}^{} P_1^T Q_F^{} F
  = F^{\ast} \widetilde{Q}_F^{\ast} \widetilde{J}^{} \widetilde{Q}_F^{} F
  = F^{\ast} J^{} F^{}.
\end{displaymath}

Since the JQR requires both row and column pivoting
(see~\cite{SingerSanja-2006}), matrix $\widetilde{G}$, with its columns
prepermuted according to $P_2$, the column pivoting of the JQR, will
then be factored by the ordinary (tall-and-skinny) QR factorization
(e.g., by the LAPACK routine \texttt{ZGEQR}).  The latter QR
factorization does not employ column pivoting, but in principle the
row pivoting or presorting may be used:
\begin{displaymath}
  P_3^{} (\widetilde{G}^{} P_2^{}) = Q_G^{} G^{}; \quad
  \widetilde{Q}_G^{\ast} \widetilde{Q}_G^{} = I_n^{}, \quad
  \widetilde{Q}_G^{} \assgn P_3^T Q_G^{},
\end{displaymath}
where $G^{} \in \mathbb{C}^{n \times n}$ is upper triangular, and
$Q_G^{} \in \mathbb{C}^{m \times n}$, which is not needed in our
application.  We also do not depend on the special forms of $F$ and
$G$ later on.

From (\ref{3.1}) it follows that the $\widetilde{J}$-GHSVD of
$\widetilde{F}$ and $\widetilde{G}$, and that of $\widetilde{F} P_2$
and $\widetilde{G} P_2$, differ only in the column permutation of the
right singular vectors, i.e., $\widetilde{X}^{} = X^{} P_2^T$, or,
from (\ref{3.2}), the row permutation of the eigenvectors, i.e.,
$\widetilde{Z} = P_2 Z$, while $\Sigma$, and thus $\Lambda$, stay the
same.  Therefore, the square factors $F$, $G$, and $J$ can be used in
place of $\widetilde{F}$, $\widetilde{G}$, and $\widetilde{J}$
throughout the rest of the computation, and then the results could be
easily converted back to the ones of the original problem.

Whenever $\widetilde{G}$ might be badly conditioned, Phase~2 could be
skipped, or we should resort to a slower but more stable QR
factorization of $\widetilde{G} P_2^{}$ with the column pivoting,
\begin{displaymath}
  P_4^{} (\widetilde{G}^{} P_2^{}) P_5^{} = Q_{G'}^{} G'; \quad
  \widetilde{Q}_{G'}^{\ast} \widetilde{Q}_{G'}^{} = I_n^{}, \quad
  \widetilde{Q}_{G'}^{} \assgn P_4^T Q_{G'}^{},
\end{displaymath}
where $P_5$ has to be applied back to $F$, and $P_4$ comes from an
optional row pivoting.  The column-pivoted QR factorization is
provided by the LAPACK routine \texttt{ZGEQP3}.

Then, $F' \assgn F^{} P_5^{}$ and $G'$ (provided that it is not
rank-deficient according to a user-defined tolerance) could be
substituted for $F$ and $G$ in the rest of the computation.  For $X$
(or $Z$) thus obtained it holds $\widetilde{X}^{} P' = X^{}$ (or,
$P' Z^{} = \widetilde{Z}^{}$), where $P' \assgn P_2^{} P_5^{}$.  Such
an approach is not required for our datasets, and therefore it was not
tested.
\subsection{$\widetilde{J}$-dot products and norms}\label{ss:5.1}
Since $\widetilde{J}$ can in principle contain the positive and the
negative signs in any order, and vectorization is strongly desired,
a $\widetilde{J}$-dot product of two vectors,
$f^{\ast} \widetilde{J}^{} g^{}$, is computed as $\mathsf{v}$
piecewise sums $\Sigma_j$, $1 \le j \le \mathsf{v}$,
\begin{align*}
  \Real(\Sigma_j) & = \Real(\Sigma_j) + \widetilde{J}_i (\Real(f_i)\Real(g_i) + \Imag(f_i)\Imag(g_i)),\\
  \Imag(\Sigma_j) & = \Imag(\Sigma_j) + \widetilde{J}_i (\Real(f_i)\Imag(g_i) - \Imag(f_i)\Real(g_i)),
\end{align*}
where $i$ starts with a value of $j$ and increments in steps of
$\mathsf{v}$ up to $m$.  The components $\Real(\Sigma)$ and
$\Imag(\Sigma)$ of the resulting $\Sigma$ are obtained by
\texttt{SUM}-reducing $\Real(\Sigma_j)$ and $\Imag(\Sigma_j)$,
respectively.  Similarly, the square of the $\widetilde{J}$-norm,
$f^{\ast} \widetilde{J}^{} f^{}$ is computed by \texttt{SUM}-reducing
$\Sigma_j'$, where
\begin{displaymath}
  \Sigma_j' = \Sigma_j' + \widetilde{J}_i^{} (\Real(f_i)^2 + \Imag(f_i)^2).
\end{displaymath}

The square of the $\widetilde{J}$-``norm'' of a vector thus obtained
can be positive or negative, with a possibility of cancellations
inadvertently occurring in the summations.  It is an open question how
to compute the squares of the $\widetilde{J}$-``norms'' both
efficiently and accurately, though one possible speed improvement
might be to encode $\widetilde{J}$ as described in
subsection~\ref{sss:4.5.1} and simplify the above three piecewise
summations accordingly.
\subsection{Pivoting}\label{ss:5.2}
\begin{figure}[hbt]
  \begin{center}
    \includegraphics{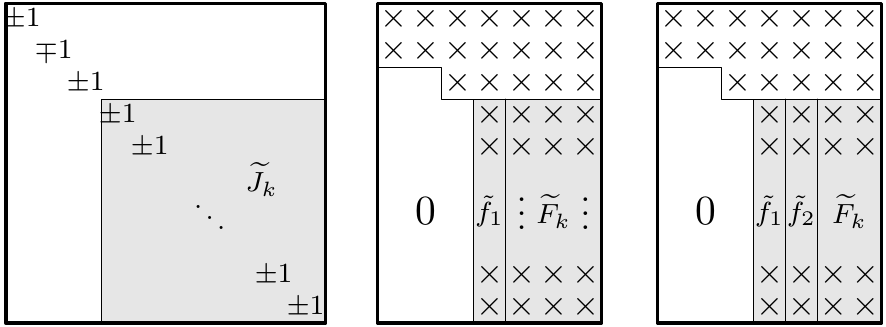}
  \end{center}
  \caption{Choosing a single pivot column or two pivot columns.
    Matrices $\widetilde{J}_k$ and $\widetilde{F}_k$ are shaded.}
  \label{fig:1}
\end{figure}

To achieve the maximal numerical stability, the JQR factorization is
usually performed with complete pivoting.  In the first step the pivot
column(s) are chosen from the $\widetilde{J}$--Grammian matrix
$H = \widetilde{F}^{\ast} \widetilde{J} \widetilde{F}$,
and later on, in the $k$th step, from the $\widetilde{J}_k$--Grammian
matrix
$H_k^{} = \widetilde{F}_k^{\ast} \widetilde{J}_k^{} \widetilde{F}_k^{}$,
where $\widetilde{F}_k$ the is a part of the matrix yet to be reduced,
and $\widetilde{J}_k^{}$ is the matrix of signs that corresponds to
the unreduced matrix $\widetilde{F}_k$ (see Fig.~\ref{fig:1}).  The
complete pivoting in the first step needs formation of the whole $H$,
i.e., $\mathcal{O}(m n^2)$ floating-point operations.  Such an
approach, consistently implemented throughout the algorithm, leads to
$\mathcal{O}(m^2 n^2)$ operations solely for the choice of pivots.
Therefore, we relaxed the pivoting strategy to the diagonal pivoting
supplemented with the partial
pivoting~\cite[Algorithm C]{Bunch-Kaufman-77}.
\subsubsection{Diagonal and partial pivoting}\label{sss:5.2.1}
First, $n - k + 1$ squares of the $\widetilde{J}_k$-norms
$h_{ii}^{[k]} \assgn \tilde{f}_i^{\ast} \widetilde{J}_k^{} \tilde{f}_i^{}$,
where $\tilde{f}_i$ is the $i$th column of $\widetilde{F}_k$ and
$i \ge 1$, are computed in a parallel-do loop over $i$ and stored in a
work array.  Since all columns have the length of $m - k + 1$, each
parallel loop iteration executes (sequentially) in approximately the
same time and the work is therefore well balanced among the threads.

Let $j \ge 1$ be the smallest index such that
$|h_{jj}^{[k]}| \ge |h_{ii}^{[k]}|$ for all $i$.  If $j > 1$, the
$k$th and the $(j + k - 1)$th columns of $\widetilde{F}$ (and thus
also the first and the $j$th columns of $\widetilde{F}_k$) are
\emph{swapped\/}.  If $k = n$, the column pivoting is completed.

Otherwise, $n - k$ $\widetilde{J}_k$-dot products
$h_{1j}^{[k]} \assgn \tilde{f}_1^{\ast} \widetilde{J}_k^{} \tilde{f}_j^{}$,
where $\tilde{f}_j$ is the $j$th column of $\widetilde{F}_k$ and
$j > 1$, are computed in a parallel-do loop over $j$ and stored in a
complex workspace, while their magnitudes $|h_{1j}^{[k]}|$ are placed
in a real workspace.  Same as above, this work is well balanced among
the threads.

Let $i > 1$ be the smallest index such that
$|h_{1i}^{[k]}| \ge |h_{1j}^{[k]}|$ for all $j > 1$.  As
in~\cite{Bunch-Kaufman-77}, if
$|h_{11}^{[k]}| \ge \alpha |h_{1i}^{[k]}|$, with
$\alpha \assgn (1 + \sqrt{17}) / 8$, the column pivoting in the step
$k$ is completed.

Otherwise, $n - k$ $\widetilde{J}_k$-dot products
$h_{il}^{[k]} \assgn \tilde{f}_i^{\ast} \widetilde{J}_k^{} \tilde{f}_l^{}$,
where $\tilde{f}_l$ is the $l$th column of $\widetilde{F}_k$ and
$i \ne l \ge 1$, are computed in a parallel-do loop over $l$ and their
magnitudes $|h_{il}^{[k]}|$ are stored in a real workspace.  This work
is only slightly imbalanced among threads, since for $l = i$ a thread
assigned to the iteration sets $|h_{ll}^{[k]}| = 0$, excluding the
value and its index from the search for a maximum, unless all other
values are also 0.

Let $j \ge 1$ be the smallest index such that
$|h_{ij}^{[k]}| \ge |h_{il}^{[k]}|$ for all $l$.  As
in~\cite{Bunch-Kaufman-77}, if
$|h_{11}^{[k]}|^{} |h_{ij}^{[k]}|^{} \ge \alpha |h_{1i}^{[k]}|^2$,
the column pivoting for the step $k$ is completed; else, if
$|h_{ii}^{[k]}| \ge \alpha |h_{ij}^{[k]}|$, the $k$th and the
$(i + k - 1)$th columns of $\widetilde{F}$ (and thus also the first
and the $i$th columns of $\widetilde{F}_k$) are \emph{swapped\/} and
the column pivoting in the step $k$ is completed.

Otherwise, a $2 \times 2$ pivot is chosen, by taking the first column
of $\widetilde{F}_k$ and \emph{swapping\/} the $(k + 1)$th and
$(i + k - 1)$th columns of $\widetilde{F}$ (and thus also the second
and the $i$th columns of $\widetilde{F}_k$), if $i \ne 2$; else, the
second pivot column is already in place.

The pivot column(s) have thus been brought to the front of the matrix
$\widetilde{F}_k$ by at most two column swaps.  The ensuing row
pivoting is explained further below.
\subsection{Hyperbolic Householder reflectors}\label{ss:5.3}
If a single pivot is chosen, the first column $\tilde{f}_1$ of
$\widetilde{F}_k$ is reduced by a hyperbolic Householder
reflector~\cite[Theorem~4.4]{SingerSanja-SingerSasa-2008} to a vector
$f_1 = c_1 e_1$, where $c_1 \in \mathbb{C}$ and $e_1$ is the first
vector of the canonical base.  A variant of
\cite[Theorem~4.4]{SingerSanja-SingerSasa-2008} for the hyperbolic
scalar product and a simple shape of $f_1$ follows.
\begin{theorem}
  Let $\widetilde{J}_k$ be a hyperbolic scalar product matrix of order
  $\ell$.  Let $\tilde{f}_1, f_1 \in \mathbb{C}^{\ell}$ be two
  distinct vectors.  There exists a basic $\widetilde{J}_k$ reflector
  $H(w)$,
  \begin{equation}
    H(w) = I - 2 w (w^{\ast} \widetilde{J}_k w)^{+} w^{\ast} \widetilde{J}_k,
    \label{5.2}
  \end{equation}
  such that $H(w) \tilde{f}_1 = f_1$ if and only if $\tilde{f}_1$ and
  $f_1$ satisfy the $\widetilde{J}_k$-isometry and
  $\widetilde{J}_k$-symmetry property, respectively
  \begin{align}
    \tilde{f}_1^{\ast} \widetilde{J}_k^{} \tilde{f}_1^{}
    & = f_1^{\ast} \widetilde{J}_k f_1^{},
    \label{5.3} \\
    \tilde{f}_1^{\ast} \widetilde{J}_k^{} f_1^{}
    & = f_1^{\ast} \widetilde{J}_k \tilde{f}_1^{},
    \label{5.4}
  \end{align}
  and $d = \tilde{f}_1 - f_1 \neq 0$ is nondegenerate, i.e.,
  $d^{\ast} \widetilde{J}_k d \neq 0$.
  Furthermore, whenever $H(w)$ exists, it is unique.  $H(w)$ can be
  generated by any $w \in \mathbb{C}^{\ell}$ such that $w = \lambda d$,
  $\lambda \in \mathbb{C}\setminus\{0\}$.
  Finally, the same remains valid if we replace $f_1$ by $-f_1$, and
  $d$ by $s = f_1 + \tilde{f}_1$.
  \label{thm5.1}
\end{theorem}

If $\tilde{f_1} = f_1$, there is nothing to do in this step, so we
take $H(w) = I$.  Otherwise, since we want to obtain $f_1$ in the form
$f_1 = c_1 e_1$, equation (\ref{5.3}) is equivalent to a requirement
that the sign of $\tilde{\jmath}_{11}^{}$, the first diagonal element
of $\widetilde{J}_k$, is equal to the sign of
$\tilde{f}_1^{\ast} \widetilde{J}_k^{} \tilde{f}_1^{}
= |c_1^{}|^2 \tilde{\jmath}_{11}^{}
= f_1^{\ast} \widetilde{J}_k^{} f_1^{}$.
As $\tilde{f}_1^{\ast} \widetilde{J}_k^{} \tilde{f}_1^{}$ has already
been computed by the diagonal pivoting (see
subsection~\ref{sss:5.2.1}), it is trivial to check if the requirement
holds.

If not, it can be shown that there exists at least one element with
the correct sign in $\widetilde{J}_k^{}$.  When there is more than one
such element, our implementation sequentially finds the one that
corresponds to the largest element in $\tilde{f}_1$ by magnitude (say,
$l$th).  By permuting the diagonal of $\widetilde{J}_k^{}$, this
element can be brought to the first diagonal position.  This implies a
corresponding row permutation of $\widetilde{F}_k$ that swaps the
first and the $l$th rows of $\widetilde{F}_k$.  For that, we employ the
parallel \texttt{ZSWAP} routine, but should the rows be short enough,
a sequential version of the routine could be considered instead.

Relation (\ref{5.3}) implies that
$|c_1| = |\tilde{f}_1^{\ast} \widetilde{J}_k^{} \tilde{f}_1^{}|^{1/2}$.
In general, $c_1$ is a complex number, $c_1 = r e^{\mathrm{i} \delta}$,
and we have already determined $r=|c_1|$.  It remains to find
$\delta = \arg(c_1)$.

From (\ref{5.4}) it follows
\begin{equation}
  \bar{\tilde{f}}_{11} \tilde{\jmath}_{11} c_1 = \bar{c}_1 \tilde{\jmath}_{11} \tilde{f}_{11},
\label{5.5}
\end{equation}
where $\tilde{f}_{11} \assgn r_1 e^{\mathrm{i} \delta_1}$ is the first
element in $\tilde{f}_1$.  Relation (\ref{5.5}) can be divided by
$\tilde{\jmath}_{11}$ and written as
$r_1 r e^{\mathrm{i} (\delta - \delta_1)} = r_1 r e^{-\mathrm{i} (\delta - \delta_1)}$.
Since $r_1, r \neq 0$, we may choose $\delta = \delta_1$ for
$\arg(c_1)$.  We only need to compute $e^{\mathrm{i} \delta}$, so
$e^{\mathrm{i} \delta} = \tilde{f}_{11} / |\tilde{f}_{11}|$.  Now we
have satisfied the conditions (\ref{5.3})--(\ref{5.4}) for
construction of a reflector that maps $\tilde{f}_1$ to $c_1 e_1$, or
to $-c_1 e_1$.

We aim to compute $d$ or $s$ accurately.  First,
$w^{\ast} \widetilde{J}_k w$ is needed, where
$w = \tilde{f}_1 \pm f_1$ (with the addition for $s$ and the
subtraction for $d$).  Since $\tilde{f}_1$ and $f_1$ satisfy
(\ref{5.3})--(\ref{5.4}),
\begin{align*}
  w^{\ast} \widetilde{J}_k w & = \big(\tilde{f}_1 \pm f_1\big)^{\ast}
  \widetilde{J}_k \big(\tilde{f}_1 \pm f_1\big)
  = 2 \big(\tilde{f}_1^{\ast} \widetilde{J}_k^{} \tilde{f}_1^{}
  \pm f_1^{\ast} \widetilde{J}_k^{} \tilde{f}_1^{} \big)
  = 2 \big( \tilde{f}_1^{\ast} \widetilde{J}_k^{} \tilde{f}_1^{}
  \pm \bar{c}_1 \tilde{\jmath}_{11} \tilde{f}_{11} \big) \\
  & = 2 \left(\tilde{f}_1^{\ast} \widetilde{J}_k^{} \tilde{f}_1^{} \pm
  |\tilde{f}_1^{\ast} \widetilde{J}_k^{} \tilde{f}_1^{}|^{1/2}
  \frac{\bar{\tilde{f}}_{11}}{|\tilde{f}_{11}|} \tilde{\jmath}_{11} \tilde{f}_{11} \right)
  = 2 \big( \tilde{f}_1^{\ast} \widetilde{J}_k^{} \tilde{f}_1^{} \pm
  |\tilde{f}_1^{\ast} \widetilde{J}_k^{} \tilde{f}_1^{}|^{1/2}
  |\tilde{f}_{11}| \tilde{\jmath}_{11} \big).
\end{align*}
Since (\ref{5.3}) should be valid, we have
$\sgn(\tilde{f}_1^{\ast} \widetilde{J}_k^{} \tilde{f}_1^{}) = \tilde{\jmath}_{11}$,
and both terms in the previous relation have the same sign.  Therefore,
to avoid unnecessary cancellation, our choice is $w = s$, and $H(s)$
from (\ref{5.2}) is then equal to
\begin{displaymath}
  H(s) = I + \tau s s^{\ast} \widetilde{J}_k, \quad
  \tau = -1 / \big(\tilde{f}_1^{\ast} \widetilde{J}_k^{} \tilde{f}_1^{} +
  |\tilde{f}_1^{\ast} \widetilde{J}_k^{} \tilde{f}_1^{}|^{1/2}
  \, |\tilde{f}_{11}| \tilde{\jmath}_{11}\big).
\end{displaymath}

The update of a column $\tilde{f_j}$, $j > 1$, is performed as
$H(s) \tilde{f}_j = \tilde{f}_j + \tau s (s^{\ast} \widetilde{J}_k \tilde{f}_j)$,
which involves computing the $\widetilde{J}_k$-dot product of $s$ and
$\tilde{f}_j$, scaling it by $\tau$, and calling the \texttt{ZAXPY}
BLAS~1 routine.  We opted for a sequential \texttt{ZAXPY} version, but
it can be argued that for the extremely long columns a parallel
version would be better suited.

The column updates are mutually independent, and therefore can be
performed in a parallel-do loop over $j$, with the work being well
balanced among the threads.

In our application it is not required that the reflector generators
($s_k$ and $\tau_k$ in the step $k$) are preserved, but that is
nevertheless done with time overhead close to none in a separate
complex matrix (the vector $s_k$ is stored in its $k$th column, with
the leading rows set to 0) and a real vector (the scalar $\tau_k$ is
stored as its $k$th element).

The situation is more complicated if a pair of pivot columns is
chosen.  According to~\cite{SingerSanja-SingerSasa-2008}, the
reduction can be performed by the block Householder matrix defined by
these two columns.  However, the computation of a block Householder
reflector is a more difficult approach than the computation of a
variant of the URV factorization with $U$ hyperbolic, not unitary.
The proposed factorization is similar, but not equal to the hyperbolic
(sometimes also called signed) URV factorization presented
in~\cite{Zhou-VanderVeen-2012}.
\begin{definition}
  Let $F \in \mathbb{C}^{m \times n}$ and
  $J \in \mathbb{Z}^{m \times m}$, $J = \diag(\pm 1)$, be given
  matrices.  Let $J' = P^T J P$ for any permutation matrix $P$.
  A factorization $F = U R V$, where $U \in \mathbb{C}^{m \times m}$
  is $J'$-unitary, $V \in \mathbb{C}^{n \times n}$ is unitary, and
  $R = \left[\begin{smallmatrix}
      R_0 \\
      0
  \end{smallmatrix}\right] \in \mathbb{C}^{m \times n}$, with
  $R_0 \in \mathbb{C}^{n \times n}$ upper triangular, is called
  a $(J, I)$ URV factorization according to $J$.
\end{definition}
Note that the $(J, I)$ URV factorization is not unique.

The Grammian matrix $A_{12}$ of the pivot block
$\widetilde{F}_{12} \assgn [\tilde{f}_1\ \tilde{f}_2]$, i.e.,
$A_{12} \assgn \widetilde{F}_{12}^{\ast} \widetilde{J}_k^{} \widetilde{F}_{12}^{}$
(see Fig.\ \ref{fig:1}) is nonsingular.  Since the pivot strategy has
chosen this block for the pivot block, the off-diagonal elements are
larger in magnitude than the diagonal elements of $A_{12}$.  The
matrix $A_{12}$ will be diagonalized by a single Jacobi rotation
$R_k$,
\begin{displaymath}
  R_k^{\ast} {A}_{12}^{} R_k^{} = R_k^{\ast} \widetilde{F}_{12}^{\ast} \widetilde{J}_k^{} \widetilde{F}_{12}^{} R_k^{}
  = \diag(d_k, d_{k+1}).
\end{displaymath}
This shows that the columns
$[\hat{f}_1\ \hat{f}_2] \assgn \widetilde{F}_{12} R_k$
are mutually $\widetilde{J}_k$ orthogonal -- their scalar product is
zero, and the squares of their $\widetilde{J}_k$ norms are $d_k$ and
$d_{k+1}$, respectively.  Since these columns are mutually
$\widetilde{J}_k$ orthogonal, the only way to handle them is applying
two successive hyperbolic Householder transformations to reduce the
pivot matrix $\widetilde{F}_{12}$ to an upper triangular matrix
$\widehat{F}_{12}$.

The first hyperbolic transformation, by $H(s_k)$, reduces the first
column to a single element.  The second hyperbolic transformation, by
$H(s_{k+1})$, then reduces the elements, already transformed by
$H(s_k)$, of the shortened (the first row of $\widehat{F}_{12}$ is not
changed anymore) second column.  We should then multiply
$\widehat{F}_{12}$ by $R_k^{\ast}$ to obtain the reduced matrix
$F_{12}$ with, generally full, topmost $2 \times 2$ block
\begin{displaymath}
  F_{12}^{} \assgn \widehat{F}_{12}^{} R_k^{\ast}
  = \begin{bmatrix}
    f_{11} & f_{21} & 0 & \cdots & 0 \\
    f_{12} & f_{22} & 0 & \cdots & 0
   \end{bmatrix}^T.
\end{displaymath}
Finally, the step counter $k$ is incremented by $2$ (instead of $1$
for a single pivot), and the process continues again with the column
pivoting, as in subsection~\ref{ss:5.2}.

When present, the multiplications from right, first by the Jacobi
rotation $R_k^{}$, and then by $R_k^{\ast}$, cancel each other (since
$R_k^{} R_k^{\ast} = I_2$), so the right matrix $V$ in this URV-like
factorization is in fact identity.

The Jacobi rotations are applied by calling the parallel
\texttt{ZROT} routine, but should the columns be short enough, a
sequential version of the routine could be warranted.

After the final step of the reduction (for $k = n$), the new $F$ and
$J$ are square matrices of order $n$.  More precisely, $J$ is the
leading part of $\widetilde{J}$ as left after all permutations due to
the row pivoting, and $F$, unlike the standard compact representation
of the LAPACK's QR factorizations, has zeros below the diagonal block
set explicitly.
\subsection{Testing}\label{ss:5.4}
In Table~\ref{tbl:2} the wall execution times of both the JQR and the
TSQR are shown.  It is evident that there is still room for the future
JQR's efficiency improvement, which might be achieved by blocking and
delaying the columns updates.
\begin{table}[hbt]
  {\footnotesize\caption{The wall execution time (wtime) in seconds of
      the Phase~2 steps: the hyperbolic QR (JQR), the tall-and-skinny
      QR (\texttt{ZGEQR} from LAPACK), and the prepermutation of the columns of
      $\widetilde{G}$ by $P_2$ (max.\ wtime with 32, 64, and \emph{48}
      threads).  In the last column is the number of $2 \times 2$
      pivots chosen.}
    \label{tbl:2}
    \begin{center}
      \begin{tabular}{@{}ccccccccc@{}}
        \toprule
        \multirow[c]{2}{*}{ID} & \multicolumn{3}{c}{JQR wtime [s]} & \multicolumn{3}{c}{TSQR wtime [s]} & $\widetilde{G}P_2$ & $2 \times 2$\\
        & 32 thr. & 64 thr. & \emph{48 thr\/}. & 32 thr. & 64 thr. & \emph{48 thr\/}. & wtime [s] & pivots\\
        \toprule
        A1 & $\hphantom{0}162.52$ & $\hphantom{0}179.02$ & $\hphantom{0}66.52$ & $10.53$ & $\hphantom{0}7.19$ & $\hphantom{0}1.80$ & $0.31$ & $11$\\
        \midrule
        A2 & $\hphantom{0}449.13$ & $\hphantom{0}482.63$ &            $184.46$ & $14.06$ & $\hphantom{0}8.26$ & $\hphantom{0}4.07$ & $0.87$ & $\hphantom{0}9$\\
        \midrule
        A3 &            $1041.92$ &            $1104.48$ &            $383.60$ & $22.01$ &            $14.71$ & $\hphantom{0}8.03$ & $1.14$ & $\hphantom{0}7$\\
        \midrule
        A4 &            $2137.63$ &            $2249.65$ &            $753.71$ & $42.82$ &            $27.38$ &            $16.54$ & $2.37$ & $\hphantom{0}6$\\
        \midrule
        B1 & $\hphantom{0}181.13$ & $\hphantom{0}177.44$ & $\hphantom{0}64.72$ & $12.00$ &            $12.21$ & $\hphantom{0}3.41$ & $0.06$ & $16$\\
        \midrule
        B2 & $\hphantom{0}520.15$ & $\hphantom{0}495.70$ &            $186.48$ & $15.63$ &            $10.29$ & $\hphantom{0}4.28$ & $0.69$ & $16$\\
        \midrule
        B3 &            $1263.68$ &            $1182.83$ &            $420.02$ & $36.66$ &            $20.22$ &            $10.29$ & $1.59$ & $20$\\
        \midrule
        B4 &            $2780.19$ &            $2439.50$ &            $841.83$ & $49.04$ &            $33.12$ &            $19.35$ & $3.41$ & $16$\\
        \bottomrule
      \end{tabular}
  \end{center}}
\end{table}
\subsubsection{Prepermuting of $\widetilde{G}$}\label{sss:5.4.1}
Table~\ref{tbl:2} also contains the wall time of preparing
$\widetilde{G} P_2$ in parallel.  The fastest way to prepermute the
columns of $\widetilde{G}$ is to copy them to another matrix of the
same size, with the column $j$ going to $\pi_2(j)$, where $\pi_2$
denotes the permutation represented by $P_2$.  Such copying occurs in
a parallel-do loop over $j$.
%
%
\section{Phase~3 -- generalized hyperbolic SVD}\label{s:6}
%
%
Vjeran Hari in his PhD thesis~\cite{Hari-84} developed a method for
solving the generalized eigenproblem, when at least one of the two
matrices is positive definite.  The method is based on the ideas from
the PhD thesis of Katharina Zimmermann~\cite{Zimmermann-69}, and has
been revisited recently in~\cite{Hari-2018,Hari-2019}.

Based on the Hari--Zimmermann algorithm for the generalized
eigenproblem, in~\cite{Novakovic-SingerSanja-SingerSasa-2015} a
one-sided method for computing the real generalized SVD has been
derived.  The main trick, how to obtain a one-sided method for the SVD
from the two-sided method for the eigenproblem is always the same:
think about the transformations in the two-sided fashion, and apply
them from one (right or left) side on a matrix factor.

Since the elements of the pivot submatrices $\widehat{H}$ of $H$ and
$\widehat{S}$ of $S$ are the scalar products of the columns of $F$ and
$G$, respectively, it is easier to write the transformations in terms
of the elements of $\widehat{H}$ and $\widehat{S}$,
\begin{equation}
  \widehat{H} = \begin{bmatrix}
    h_{pp} & h_{pq} \\
    \bar{h}_{pq} & h_{qq}
  \end{bmatrix}
  = \begin{bmatrix}
    f_p^{\ast} J f_p^{} & f_p^{\ast} J f_q^{} \\
    \overline{f_p^{\ast} J f_q^{}} & f_q^{\ast} J f_q^{}
  \end{bmatrix}, \quad
  \widehat{S} = \begin{bmatrix}
    s_{pp} & s_{pq} \\
    \bar{s}_{pq} & s_{qq}
  \end{bmatrix}
  = \begin{bmatrix}
    g_p^{\ast} g_p^{} & g_p^{\ast} g_q^{} \\
    \overline{g_p^{\ast} g_q^{}} & g_q^{\ast} g_q^{}
  \end{bmatrix},
  \label{6.1}
\end{equation}
instead of in terms of the columns $f_p^{}$, $f_q^{}$, $g_p^{}$, and
$g_q^{}$.

The original method consists of $3$ active transformations, and an
auxiliary transformation that helps in coupling them all together.
\subsection{Pointwise algorithm}\label{ss:6.1}
The whole pointwise algorithm (with a pair of $2 \times 2$ pivot
submatrices in each annihilation step) is taken from the PhD thesis of
Vjeran Hari~\cite{Hari-84}.  However, we feel that the algorithm
should be presented succinctly here, to aid its implementors, and also
to incorporate some minor corrections.
\subsubsection{Preprocessing}\label{sss:6.1.1}
In the preprocessing step, $H$ and $S$ are scaled by a diagonal matrix
$D$ such that $\diag(DSD) = I$, i.e.,
\begin{displaymath}
  H_0^{} \assgn D H D, \quad S_0^{} \assgn D S D, \quad
  D = \diag\big(s_{11}^{-1/2}, s_{22}^{-1/2}, \ldots, s_{nn}^{-1/2}\big).
\end{displaymath}

Such preprocessing can be done only once, at the start of the
algorithm, or it can be done before each annihilation step for the
pivot column pair in question.  The latter might seem redundant, but
in a floating-point realization of the algorithm, after enough steps,
$\diag(DSD)$ could veer off the identity enough to warrant such a
rescaling.  In that case, form a matrix $D_0$ that has, for the chosen
pivot indices $(p,q)$, as its $p$th and $q$th diagonal entries
$(s_{pp})^{-1/2}$ and $(s_{qq})^{-1/2}$, respectively, while being
equal to the identity elsewhere.  For the approach with a single
prescaling, let $D_0 = I$.

In both cases, let $\widehat{D}_0$ be a $2 \times 2$ restriction of
$D_0$ to the $p$th and the $q$th rows and columns.  We have
implemented the pivot pair prescaling in each annihilation step.
\subsubsection{Diagonalization of $\widehat{S}_0$}\label{sss:6.1.2}
In the first step the $2 \times 2$ pivot submatrix $\widehat{S}_0^{}$
of $S_0^{}$ (at the crossings of the $p$th and the $q$th rows and
columns) is diagonalized by a complex Jacobi rotation $\widehat{R}_1$,
where $\widehat{R}_k$ for $k\ge 1$ is
\begin{equation}
  \widehat{R}_k = \begin{bmatrix}
    \hphantom{-}\cos \varphi_k & e^{\mathrm{i} \alpha_k} \sin \varphi_k\\
    -e^{-\mathrm{i} \alpha_k} \sin \varphi_k & \cos \varphi_k
  \end{bmatrix}.
  \label{6.2}
\end{equation}
The same transformation is then applied to $H_0$, to keep the new pair
equivalent to the original one.  After that, the new pair is
$(H_1^{}, S_1^{}) \assgn (R_1^{\ast} H_0^{} R_1^{}, R_1^{\ast} S_0^{} R_1^{})$,
where $R_1^{} = I$, except at the pivot positions, where
$R_1^{} = \widehat{R}_1^{}$.  Since $H$ and $S$ have been preprocessed
as in subsection~\ref{sss:6.1.1}, the diagonal elements of
$\widehat{S}_0$ are the same, and we may choose $\varphi_1 = -\pi/4$
in (\ref{6.2}).  Now it is easy to determine that
$\alpha_1 = \arg(s_{pq})$.  If $s_{pq}$, written in the trigonometric
form, was $s_{pq} = x e^{\mathrm{i} \alpha_1}$, with $x = |s_{pq}|$,
before the transformation, then after it we obtain
\begin{equation}
  \widehat{S}_1^{} = \diag(1 + x, 1 - x).
  \label{6.3}
\end{equation}
\subsubsection{Rescaling of $\widehat{S}_1$}\label{sss:6.1.3}
The second step rescales the diagonal of $S_1^{}$ to ones, and
rescales $H_1^{}$ with the same diagonal matrix.  After the
transformation, similar to the preprocessing step, we obtain
$H_2^{} \assgn D_2^{} H_1^{} D_2^{}$, $S_2^{} \assgn D_2^{} S_1^{} D_2^{}$.
From (\ref{6.3}) we conclude that
$\widehat{D}_2^{} = \diag ( (1 + x)^{-1/2}, (1 - x)^{-1/2})$.
Elements of $\widehat{H}_2^{}$ are
\begin{equation}
  \widehat{H}_2^{} = \begin{bmatrix}
    h_{pp}^{(2)} & h_{pq}^{(2)} \\
    \bar{h}_{pq}^{(2)} & h_{qq}^{(2)}
  \end{bmatrix}
  = \frac{1}{2} \begin{bmatrix}
    \frac{1}{1 + x} (h_{pp} + h_{qq} + 2u)
    & \frac{e^{\mathrm{i} \alpha_1}}{\sqrt{1 - x^2}} (h_{qq} - h_{pp} + 2 \mathrm{i} v) \\[3pt]
    \frac{e^{-\mathrm{i} \alpha_1}}{\sqrt{1 - x^2}}(h_{qq} - h_{pp} - 2 \mathrm{i} v)
    & \frac{1}{1 - x} (h_{pp} + h_{qq} - 2u)
  \end{bmatrix},
  \label{6.4}
\end{equation}
where
\begin{equation}
  u + \mathrm{i}v = e^{-\mathrm{i} \alpha_1} h_{pq} = e^{-\mathrm{i} \arg(s_{pq})} h_{pq}.
  \label{6.5}
\end{equation}
\subsubsection{Diagonalization of $\widehat{H}_2$}\label{sss:6.1.4}
In the third step the pivot submatrix $\widehat{H}_2^{}$ of $H_2^{}$
is diagonalized by a complex Jacobi rotation $\widehat{R}_3$ of the
form~(\ref{6.2}).  The third transformation can be written as
$H_3^{} = R_3^{\ast} H_2^{} R_3^{}$, $S_3^{} = R_3^{\ast} S_2^{} R_3^{}$,
where $R_3^{} = I$, except at the pivot positions, where
$R_3^{} = \widehat{R}_3^{}$.  Then $\varphi_3$ in~(\ref{6.2}) is
written as $\varphi_3 = \vartheta + \pi/4$ to express the
transformations in terms of $\vartheta$.  The relations for the angles
of a rotation that annihilates the off-diagonal element of
$\widehat{H}_2$ are given~by
\begin{equation}
  \tan \left( 2 \vartheta + \frac{\pi}{2} \right)
  = \sigma \frac{2 |h_{pq}^{(2)}|}{h_{qq}^{(2)} - h_{pp}^{(2)}}, \qquad
  \alpha_3 = \alpha_1 + \arg \left( \frac{h}{2} + \mathrm{i} v \right)
  + (1 - \sigma) \frac{\pi}{2},
\label{6.6}
\end{equation}
where $h = h_{qq} - h_{pp}$.  The requirement
$\gamma \assgn \alpha_3 - \alpha_1 \in \left(-\pi/2, \pi/2\right]$
yields $\sigma = \sgn(h)$.  This choice of $\sigma$ as the sign of $h$
and the constraint $\vartheta \in \left(-\pi/4, \pi/4\right]$ ensures
the convergence of the algorithm (see~\cite{Hari-84}).

Since
$\tan(2\vartheta + \pi/2) = -\cot(2\vartheta) = -1/\tan(2\vartheta)$,
from (\ref{6.4}) and (\ref{6.6}), we obtain
\begin{equation}
  \tan(2\vartheta) = \sigma \frac{2u - (h_{pp} + h_{qq}) x}
    {t \sqrt{h^2 + 4v^2}}, \qquad t \assgn \sqrt{1 - x^2}.
    \label{6.7}
\end{equation}

After the first three steps, the pivot submatrix $\widehat{S}_3^{}$ is
still diagonal (in fact identity), i.e.,
$\widehat{S}_3^{} = \widehat{R}_{}^{\ast} \widehat{S}_0^{} \widehat{R}^{} = I^{}$,
where
\begin{equation}
  \widehat{R} = \widehat{R}_1 \widehat{D}_2 \widehat{R}_3.
\label{6.8}
\end{equation}
This constructively shows that $\widehat{R}$ diagonalizes the matrix
pair $(\widehat{H}_0, \widehat{S}_0)$.
\subsubsection{Forming $\widehat{R}$ and $\widehat{Z}'$}\label{sss:6.1.5}
Hari in~\cite[Theorem~2.2]{Hari-84} has proved the form of the general
nonsingular matrix that diagonalizes a $2 \times 2$ Hermitian positive
definite matrix.  The intention of representing $\widehat{R}$ as
\begin{equation}
  \widehat{R} = \frac{1}{t}
  \begin{bmatrix}
    \hphantom{-}\cos \varphi & e^{\mathrm{i} \alpha} \sin \varphi\\
    -e^{-\mathrm{i} \beta} \sin \psi & \cos \psi
  \end{bmatrix}
  \diag(e^{\mathrm{i} \sigma_p}, e^{\mathrm{i} \sigma_q})
  \label{6.9}
\end{equation}
is to simplify (\ref{6.8}).  Comparing the elements of (\ref{6.8}) and
(\ref{6.9}) we obtain
\begin{equation}
  \begin{aligned}
  \cos \varphi & = \big(1/\sqrt{2}\big) \cdot \sqrt{1 + x \sin(2 \vartheta) +
    t \cos \gamma \cos(2 \vartheta)}, \quad 0 \leq \varphi < \pi/2,\\
  \cos \psi & = \big(1/\sqrt{2}\big) \cdot \sqrt{1 - x \sin(2 \vartheta) + t
    \cos \gamma \cos(2 \vartheta)}, \quad 0 \leq \psi < \pi/2,\\
  e^{\mathrm{i} \alpha} \sin \varphi &
  = e^{\mathrm{i} \alpha_1} \cdot ((\sin(2\vartheta) - x)
    + \mathrm{i} t \sin \gamma \cos(2\vartheta))\,/\,(2 \cos \psi),\\
  e^{-\mathrm{i} \beta} \sin \psi &
  = e^{-\mathrm{i} \alpha_1} \cdot ((\sin(2 \vartheta) + x)
    - \mathrm{i} t \sin \gamma \cos(2 \vartheta))\,/\,(2 \cos \varphi).
  \end{aligned}
  \label{6.10}
\end{equation}

Since in (\ref{6.10}) we need only $\sin \gamma$ and $\cos \gamma$,
from (\ref{6.6}) it then follows
\begin{equation}
  \tan \gamma = 2 v/h.
  \label{6.11}
\end{equation}

The fourth step only deals with a formal simplification of
$\widehat{R}$, by introducing a transformation $\Phi_4$ such that
$H_4^{} = \Phi_4^{\ast} H_3^{} \Phi_4^{}$,
$S_4^{} = \Phi_4^{\ast} S_3^{} \Phi_4^{}$.
The matrix $\Phi_4$ is a diagonal matrix equal to identity, except at
pivot positions, where
$\widehat{\Phi}_4 = \diag(e^{-\mathrm{i} \sigma_p}, e^{-\mathrm{i} \sigma_q})$.
Obviously, if $\widehat{R}$ diagonalizes the pair
$(\widehat{H}_0, \widehat{S}_0)$, then the transformation
$\widehat{Z} = \widehat{R} \widehat{\Phi}_4$ will leave the final
diagonal matrices intact.  Then, let
$\widehat{Z}' = \widehat{D}_0 \widehat{Z}$.
\subsubsection{Exceptional cases}\label{sss:6.1.6}
There can be a few exceptions in the computations of the elements of
the matrix $\widehat{Z}$ which have to be accounted for in the
algorithm.

If $h_{pq} = s_{pq} = 0$, we set $\widehat{Z} = I$, since both pivot
submatrices are already diagonal.  We could still apply the scaling by
$\widehat{D}_0$, but not count that as a transformation.

If $h_{pq} \neq 0$, but $s_{pq} = 0$ (i.e., $x = 0$), we set
$\alpha_1 = 0$ and proceed as described above to determine
$\widehat{Z}$ as an ordinary Jacobi rotation that diagonalizes
$\widehat{H}_0$.

If $h = v = 0$ in~(\ref{6.11}), i.e., when
$\arg{s_{pq}} = \arg{h_{pq}}$ and $h_{pp} = h_{qq}$, it can be shown
that $\widehat{R}_1$, the Jacobi rotation that diagonalizes
$\widehat{S}_0$, also diagonalizes $\widehat{H}_0$, so
$\widehat{Z} = \widehat{R}_1 \widehat{D}_2$.
\subsubsection{Convergence criterion and finalization}\label{sss:6.1.7}
As in~\cite{Novakovic-SingerSanja-SingerSasa-2015}, the convergence
criterion in floating-point arithmetic has to take into account the
relative magnitudes of the off-diagonal elements, compared to the
diagonal ones, in the pivot pairs.  Therefore, a pivot pair undergoes
the transformation if, for the machine precision $\varepsilon$,
\begin{equation}
  |h_{pq}| \ge \left(\max\{|h_{pp}|, |h_{qq}|\} \cdot \left(\varepsilon \sqrt{n}\right)\right)\cdot \min\{|h_{pp}|, |h_{qq}|\}
  \quad \text{or} \quad |s_{pq}| \ge \varepsilon \sqrt{n}.
  \label{6.12}
\end{equation}

If $\widehat{Z}'$ turns out to be identity, it is not applied.  If
$\cos\varphi = \cos\psi = 1$, such transformation is considered
``small'', and ``big'' otherwise.  Near the end of the process, the
transformations turn out to be small, and (in a blocking variant, the
last level of) the algorithm is stopped when no big transformations
are encountered in a sweep, to avoid perpetually applying the
transformations that spring only from the accumulation of the rounding
errors.  With blocking, the inner level(s) of the algorithm count all
transformations (big and small) in a sweep for stopping.  For details,
see~\cite{Novakovic-2015,Novakovic-SingerSanja-SingerSasa-2015}.

\paragraph{Outputs}%
The algorithm stops when the columns of the in-place transformed $F'$
and $G'$ are numerically mutually $J$-orthogonal, and orthogonal,
respectively.  Let $Z'$ be the accumulated product of the applied
transformations.  Then,
\begin{displaymath}
  \Sigma_F' = \diag\big(|f_1^{\ast} J^{} f_1^{}|^{1/2}, \ldots, |f_{\vphantom{1}n}^{\ast} J^{} f_{\vphantom{1}n}^{}|^{1/2}\big), \quad
  \Sigma_G' = \diag\big((g_1^{\ast} g_1^{})^{1/2}, \ldots, (g_n^{\ast} g_n^{})^{1/2}\big),
\end{displaymath}
and $U^{} = F' \Sigma_F^{\prime -1}$,
$V^{} = G' \Sigma_G^{\prime -1}$,
$\Sigma_j^{} = \left((\Sigma_F')_j^2 + (\Sigma_G')_j^2\right)^{1/2}$,
$\Sigma = \diag(\Sigma_1, \ldots, \Sigma_n)$,
$\Sigma_F^{} = \Sigma_F' \Sigma^{-1}$,
$\Sigma_G^{} = \Sigma_G' \Sigma^{-1}$,
$Z^{} = Z' \Sigma^{-1}$.  For the inner levels of blocking, only the
matrix $Z$ (and therefore also $\Sigma_F'$, $\Sigma_G'$, and $\Sigma$,
but not $U$ and $V$) is required.

\paragraph{Accumulating $Z^{\prime -1}$}%
Optionally, $Z^{\prime -1}$ could be obtained by accumulating
transformations $\widehat{Z}'$ from the right and their inverses
$\widehat{Z}^{\prime -1}$ from the left.  Then,
$X^{} = Z^{-1} = \Sigma^{} $, so Phase~4 would not be
needed for the full G(H)SVD.  This has not been implemented, since our
main concern is a solution of the generalized eigenproblem.
\subsubsection{The G(H)SVD algorithm}\label{sss:6.1.8}
After obtaining the matrix $\widehat{Z}'$, the pointwise implicit
Hari--Zimmermann G(H)SVD algorithm can be written similarly
to~\cite[Algorithm~3.1]{Novakovic-SingerSanja-SingerSasa-2015}.
In Algorithm~\ref{alg:phz} we take into account the signature matrix
$J$, since $H^{}=F^{\ast} J^{} F$, but with $J=I$ it reduces to a
generalized SVD method.
\begin{algorithm}
  \caption{The pointwise implicit Hari--Zimmermann G(H)SVD algorithm.}
  \label{alg:phz}
  \begin{algorithmic}
    \STATE{$Z' = I$; optional computing of $D$ and prescaling
      $H_0 = DHD$, $S_0 = DSD$, $Z' = D$;}
    \STATE{$it = 1$;}$\,$\COMMENT{The sweep counter.  Maximal
      number of sweeps, $\mathsf{C}_{\max}$, is usually $\approx 30$.}
    \REPEAT
    \FORALL{pairs $(p, q)$, $1 \leq p < q \leq n$}
    \STATE{compute $\widehat{H}$ and $\widehat{S}$ from (\ref{6.1});}
    \STATE{compute the elements of $\widehat{Z}'$ from (\ref{6.5}), (\ref{6.7}), (\ref{6.10})--(\ref{6.11});}
    \STATE{$[f_p, f_q] = [f_p, f_q] \cdot \widehat{Z}';\quad [g_p, g_q] = [g_p, g_q] \cdot \widehat{Z}';\quad [z_p', z_q'] = [z_p', z_q'] \cdot \widehat{Z}'$;}
    \ENDFOR
    \STATE{$it = it + 1$;}
    \UNTIL{(no transformations in this sweep) \OR ($it > \mathsf{C}_{\max}$)}
    \STATE{Output: $Z$ and (optionally) $U$, $V$, $\Sigma_F$,
      $\Sigma_G$, and $\Sigma$;}
  \end{algorithmic}
\end{algorithm}

If the prescaling as described in subsection~\ref{sss:6.1.1} is
performed only once, then in Algorithm~\ref{alg:phz} only
$g_p^{\ast} g_q^{}$ is computed, since the diagonal elements of
$\widehat{S}$ are assumed to be unity.  Otherwise, by passing once
through the columns $g_p$ and $g_q$ all three dot products can be
formed.  Similar holds for $f_p$, $f_q$, and the three $J$-dot
products of $\widehat{H}$.
\subsection{Vectorization}\label{ss:6.2}
Many computational building blocks of Algorithm~\ref{alg:phz} provide
both the challenges and the opportunities for the SIMD vectorization.
The most obvious such primitives are the $J$-dot products (and norms),
which are computed combining the approaches presented in
subsections~\ref{ss:5.1} (for the unstructured patterns of signs in
$J$) and \ref{sss:4.5.1} for a compact, partitioned representation of
$J$ with only a number $n_+$ of the leading positive (and therefore
$n-n_+$ tailing negative) signs.
\subsubsection{\texttt{ZVROTM}}\label{sss:6.2.1}
The column updates by $\widehat{Z}'$ cannot be realized by a single
call to a BLAS or a LAPACK routine (e.g., \texttt{ZROT}), since there
are two angles involved in a transformation, and the columns are meant
to be transformed in-place (overwritten).  For that purpose, a
\texttt{ZVROTM} routine has been implemented as a vectorized loop,
that resembles a simplified version of the BLAS routine
\texttt{DROTM}, but with the complex sines.
\subsubsection{Transformations}\label{sss:6.2.2}
The greatest challenge lies in computing the transformations in the
SIMD-parallel way, where each vector lane $i$ computes
$\widehat{Z}_i'$ for its own pivot pair with indices $(p_i, q_i)$, as
in~\cite{Novakovic-2017} for the rotations in the Jacobi SVD method.
Assume that $\mathsf{v}$ pivot pairs $(\widehat{H}_i, \widehat{S}_i)$
have been obtained, with their corresponding $p_i$ and $q_i$ indices
all different.  If there are fewer than $\mathsf{v}$ such pairs, let
$\widehat{H}_i = \widehat{S}_i = I$ for the missing indices $i$.  A
vector (e.g., $\widehat{S}_{12}^{[i]} = g_{p_i}^{\ast} g_{q_i}^{}$)
has to be kept in an array of length $\mathsf{v}$, properly aligned in
memory, in which the $i$th position holds data for the $i$th lane.
To help the compiler, the complex arithmetic operations have been
written in terms of the real and the imaginary parts of the complex
numbers in the vectorizable regions.

For the start, check for which lanes their pivot pair has to be
transformed by evaluating the criterion~(\ref{6.12}) in each lane.  To
aid the compiler, (\ref{6.12}) can be rewritten as a branch-free
arithmetic expression that has a non-zero value if and only if the
criterion is fulfilled.  If those values constitute a zero vector, no
transformation is required for any $i$, and a fresh set of pivot pairs
(if any remain) should be considered.

Then, compute $\widehat{Z}_i'$ unconditionally (i.e., for all $i$).
The idea is that the unneeded computation comes at no cost, while its
results can be discarded afterwards.  However, the exceptional cases
from subsection~\ref{sss:6.1.6} should be carefully dealt with, since
a na\"{\i}ve branching might spoil the vectorization opportunities.
The logical conditions are therefore arithmetized, while halting on
the arithmetic exceptions is suppressed.

We assume that the second argument of the intrinsics \texttt{MIN} and
\texttt{MAX} is returned when their first argument is a \texttt{NaN}.
Such a behavior is not mandated by the Fortran standard (contrary to
\texttt{fmin} and \texttt{fmax} in C), but is checked for at runtime
by our code.

Let
$\text{\texttt{B(i)}}=\text{\texttt{(1-MAX(V(i)/V(i),0))*(1-MAX(H(i)/H(i),0))}}$
and note that $\texttt{B(i)}$ is 1 if, in~(\ref{6.11}),
$v_i = h_i = 0$ (so $\tan\gamma_i = \mathtt{NaN}$), and also if
$h_i = \pm\infty$ (due to an overflow) with $v_i = 0$; otherwise, it
is 0.  The computation resumes regardless of the value of
\texttt{B(i)}.  However, $\widehat{Z}_i'$ in those lanes where
$\text{\texttt{B(i)}} = 1$ is useless, so a new, correct
$\widehat{Z}_i'$ is taken according to the rules of
subsection~\ref{sss:6.1.6} at the end, sequentially for all such $i$.
This is the only situation when $\widehat{Z}_i'$ is not computed in
one go for all $i$, but it should occur rarely in practice.  A
branch-free exception handling is also requried when
$\widehat{S}_{12}^{[i]} = 0$, since in calculating the polar form
$\big|\widehat{S}_{12}^{[i]}\big| e^{\mathrm{i} \phi_i}$ of
$\widehat{S}_{12}^{[i]}$ the divisions of the real and the imaginary
parts by the absolute value result in \texttt{NaN}s.  If $\cos\phi_i$
is denoted by \texttt{C(i)}, then selecting a default of
$\cos\phi_i = 1$ can be done by presetting the variable to 1, and
taking $\text{\texttt{C(i)}} = \text{\texttt{MIN(x(i),C(i))}}$, where
\texttt{x} is obtained by the vector division, or the multiplication
by the reciprocal.  For $\sin\phi_i$, the default value can be also
set to 1, and the correct value of 0 is obtained by multiplying the
intermediate result by a variable, set in the process of
checking~(\ref{6.12}), that is 0 if $\widehat{S}_{12}^{[i]} = 0$,
and 1 otherwise.
\subsection{Parallelization}\label{ss:6.3}
There are $N \assgn n(n-1)/2$ pairs of indices $(p_i,q_i)$ such that
they belong to a strictly upper triangle (i.e.,
$1 \le p_i < q_i \le n$) of the square matrices of order $n$.  Let
$\mathsf{P} \assgn \{(p_i,q_i) : 1 \le i \le N\}$ be a set of those
index pairs.  At most $\left\lfloor n/2\right\rfloor$ such pairs can
be chosen from $\mathsf{P}$ so that all their indices are distinct.

Let $\mathsf{S}_j$, for some $j \ge 1$, be a set of at most
$\left\lfloor n/2\right\rfloor$ index pairs, with all indices
distinct.  Then, the pivot pairs formed from the columns of $F$ and
$G$ and indexed by the elements of $\mathsf{S}_j$ can be transformed
concurrently.  We call $\mathsf{S}_j$ the $j$th parallel Jacobi step,
and a sequence of steps
$\mathsf{S} \assgn (\mathsf{S}_1, \mathsf{S}_2, \ldots, \mathsf{S}_{\bar{n}})$
a parallel (quasi-)cyclic Jacobi strategy if
$\bigcup_j \mathsf{S}_j = \mathsf{P}$, and under assumption that the
sequence is repeated forever in principle (or, in practice, until the
convergence criteria are met).  A strategy is called cyclic if its
steps are mutually disjoint; else, it is called quasi-cyclic.  In a
cycle (also called a sweep) all pivot pairs are accessed (at least
once, but maybe more under a quasi-cyclic strategy), in
$\bar{n} \ge n - 1$ steps.  We aim for the steps as large as possible.
\subsubsection{Parallel strategies}\label{sss:6.3.1}
Two classes of the parallel Jacobi strategies were under test: the
modified modulus strategy (\textsf{MM}), described in,
e.g.,~\cite{Novakovic-SingerSanja-SingerSasa-2015}, and a
generalization of the Mantharam--Eberlein
strategy~\cite{Mantharam-Eberlein-93} (\textsf{ME}), described
in~\cite{Novakovic-2015}.  The former is a quasi-cyclic strategy with
$\bar{n} = n$, but easily generated on-the-fly as the computation
progresses.  The latter is cyclic, attains $\bar{n} = n - 1$, has
provided a faster execution than \textsf{MM} to the one-sided Jacobi
SVD, with more accurate results, but is not readily available for all
even $n$, and its convergence has not yet been proven.
\subsubsection{Bordering}\label{sss:6.3.2}
For an odd $n$, no more than $(n - 1)/2$ pairs fit into a step, so
$\bar{n} \ge n$.  We therefore consider the strategies for even $n$
only, with the steps of size $n/2$, and when necessary border the
matrices by appending a zero column and a zero row, except for the new
element at position $(n+1,n+1)$, which is set to unity.
\subsubsection{Vector-Parallel algorithm}\label{sss:6.3.3}
Let a step $\mathsf{S}_j$, with $\mathsf{k}_1 \assgn n/2$ index
pairs, be given.  Then, partition $\mathsf{S}_j$ into
$\mathsf{k}_{\mathsf{v}} \assgn \left\lceil\mathsf{k}_1/\mathsf{v}\right\rceil$
disjoint subsets $\mathsf{V}_k$, where each subset has at most
$\mathsf{v}$ pairs, i.e.,
$\mathsf{S}_j = (\mathsf{V}_1, \mathsf{V}_2, \ldots, \mathsf{V}_{\mathsf{k}_{\mathsf{v}}})$.
For each subset, the requirements for the vectorized computation of
$\widehat{Z}'$ as described in subsection~\ref{ss:6.2} are satisfied.

Now, let $\mathsf{t} \ge 1$ be a number of available (OpenMP)
threads.  Then, $\mathsf{S}_j$ is traversed with a parallel-do loop
over the subsets, where a thread takes a chunk of subsets to be
processed independently, while the subsets within a chunk are handled
in sequence.

A sweep of such Vector-Parallel (VP) variant is shown in
Algorithm~\ref{alg:vp}.
\begin{algorithm}
  \caption{A sweep of the Vector-Parallel implicit HZ algorithm.}
  \label{alg:vp}
  \begin{algorithmic}
    \FORALL[a sequential loop over the steps of $\mathsf{S}$]{steps $\mathsf{S}_j \in \mathsf{S}$, $1 \le j \le \bar{n}$}
    \FORALL[an OpenMP parallel do with $\mathsf{t}$ threads]{subsets $\mathsf{V}_k \in \mathsf{S}_j$, $1 \le k \le \mathsf{k}_{\mathsf{v}}$}
    \FORALL[a sequential loop over $\mathsf{V}_k$]{$(p_i,q_i) \in \mathsf{V}_k$, $1 \le i \le |\mathsf{V}_k| \le \mathsf{v}$}
    \STATE{get
      $\widehat{H}_i = \begin{bmatrix}
        f_{p_i}^{\ast} J f_{p_i}^{} & f_{p_i}^{\ast} J f_{q_i}^{} \\
        & f_{q_i}^{\ast} J f_{q_i}^{}
      \end{bmatrix};
      \widehat{S}_i = \begin{bmatrix}
        g_{p_i}^{\ast} g_{p_i}^{} & g_{p_i}^{\ast} g_{q_i}^{} \\
        & g_{q_i}^{\ast} g_{q_i}^{}
      \end{bmatrix};$}
    \COMMENT{vectorized $a^{\ast}(J)b$}
    \ENDFOR
    \FORALL[a SIMD parallel do over $\mathsf{V}_k$]{$(p_i,q_i) \in \mathsf{V}_k$, $1 \le i \le |\mathsf{V}_k| \le \mathsf{v}$}
    \STATE{check the transformation criterion~(\ref{6.12});}
    \ENDFOR
    \STATE{if no pivot pairs have to be transformed, \textbf{cycle};}$\ $\COMMENT{a reduction}
    \FORALL[a SIMD parallel do over $\mathsf{V}_k$]{$(p_i,q_i) \in \mathsf{V}_k$, $1 \le i \le |\mathsf{V}_k| \le \mathsf{v}$}
    \STATE{compute the elements of $\widehat{Z}_i'$ from (\ref{6.5}), (\ref{6.7}), (\ref{6.10})--(\ref{6.11});}
    \ENDFOR
    \FORALL[a sequential loop over $\mathsf{V}_k$]{$(p_i,q_i) \in \mathsf{V}_k$, $1 \le i \le |\mathsf{V}_k| \le \mathsf{v}$}
    \STATE{check if $\widehat{Z}_i'$ has to be corrected and \textbf{cycle} if $\widehat{Z}_i' = I;$}
    \STATE{$[f_{p_i}, f_{q_i}] = [f_{p_i}, f_{q_i}] \cdot \widehat{Z}_i';\quad [g_{p_i}, g_{q_i}] = [g_{p_i}, g_{q_i}] \cdot \widehat{Z}_i';\quad [z_{p_i}', z_{q_i}'] = [z_{p_i}', z_{q_i}'] \cdot \widehat{Z}_i'$;}
    \COMMENT{Two \texttt{ZVROTM}s of length $m$ and one of length $n$.}
    \ENDFOR
    \ENDFOR
    \ENDFOR
  \end{algorithmic}
\end{algorithm}
Note that each thread has to have a private set of vector variables,
which is most easily done by reserving a
$\mathsf{v} \times \mathsf{t}$ rank-2 array for each variable and
making the $l$th thread access the $l$th column of such an array,
where $l$ is a thread's unique number, $1 \le l \le \mathsf{t}$.
\subsection{Blocking}\label{ss:6.4}
To better exploit the memory hierarchy by keeping data in the cache(s)
longer, a multilevel blocking principle can be applied, with the
Level~1 being the pointwise algorithm in its VP variant.  In the next,
second level of the algorithm the block columns of width
$\mathsf{w} \ge 1$ take place of the single columns, and the
$2 \times 2$ pivot pairs are replaced by
$(2\mathsf{w}) \times (2\mathsf{w})$ block pivots.  The same principle
can be applied recursively further (e.g., see~\cite{Novakovic-2015}),
but we consider only the Level~2 algorithms here.
\subsubsection{Block-column partitioning}\label{sss:6.4.1}
There are two ways a matrix can be partitioned into block columns.
The first one is to prescribe $\mathsf{w}$, and then split the matrix
into at least $\left\lceil n/\mathsf{w}\right\rceil$ block columns of
width at most $\mathsf{w}$, bearing in mind that the number of block
columns has to be even, as explained in subsection~\ref{sss:6.3.2},
and reducing the maximal width accordingly.  For $\mathsf{w} \ge 2$,
all block columns can be made to contain either $\mathsf{w}$ or
$\mathsf{w} - 1$ (but no less) columns by redistributing their
individual widths.

The second way, and the one we have chosen to implement, is to query
at run-time a number $\mathsf{t}$ of threads to be used.  A thread $l$
is to be assigned one pair of block columns with the block indices
$(\mathsf{p}_{\mathsf{j}l},\mathsf{q}_{\mathsf{j}l})$ in the
$\mathsf{j}$th block step, so the maximal width $\mathsf{w}$ is
computed as $\left\lceil n/(2\mathsf{t})\right\rceil$, with exactly
$2\mathsf{t}$ block columns.  The block column widths
$\mathsf{w}_{\mathsf{i}}$ are non-increasing across the whole
partition and are equal to either $\mathsf{w}$ or $\mathsf{w} - 1$.

The blocking overhead can dominate the actual computation time for the
matrices small enough, so we assume that $n > 2\mathsf{t}$, and
suggest a pointwise algorithm otherwise.

Each thread allocates a contiguous storage (in its own NUMA region,
but visible to the other threads) for $2\mathsf{w}$ columns (two block
columns of the maximal width) of $F$, and similarly for $G$ and $Z$.
The same amount of memory, and of the same shape, is additionally
allocated for the columns of the ``shadow'' matrices $\mathsf{F}$,
$\mathsf{G}$, and $\mathsf{Z}$.  Let, for the $l$th thread in the
$\mathsf{j}$th step of the $\mathsf{s}$th block sweep, the contents of
that storage be named
$\mathsf{X}_{\mathsf{j}l}^{[\mathsf{s}]} \assgn
[\mathsf{X}_{\mathsf{p}_{\mathsf{j}l}}^{[\mathsf{s}]} \ \mathsf{X}_{\mathsf{q}_{\mathsf{j}l}}^{[\mathsf{s}]}]$,
with $\mathsf{X} \in \{F, G, Z, \mathsf{F}, \mathsf{G}, \mathsf{Z}\}$.
If $\mathsf{w}_{\mathsf{p}_{\mathsf{j}l}} = \mathsf{w} - 1$, the
columns of $\mathsf{X}_{\mathsf{p}_{\mathsf{j}l}}$ are placed from the
second column of $\mathsf{X}_{\mathsf{j}l}$, while the columns of
$\mathsf{X}_{\mathsf{q}_{\mathsf{j}l}}$ are always placed from the
column $\mathsf{w} + 1$ of $\mathsf{X}_{\mathsf{j}l}$.  That way all
the columns of a block column pair are stored contiguously and can be
viewed by the BLAS routines as a single matrix.

The storage for the block pivots and for the workspaces, along with a
copy of $J$ and the strategy tables, is also preallocated per thread,
in MCDRAM if possible.  There are two strategy tables; one for the
outer, Level~2 Jacobi strategy $\mathsf{S}^{[2]}$, and one for the
inner, Level~1 strategy $\mathsf{S}^{[1]}$, which do not have to
belong to the same class.  The tables are fully initialized, with all
the (block) steps, before the start of the iterations.

The initial data for $F_{1l}^{[1]}$ and $G_{1l}^{[1]}$ is loaded from
the block columns
$(\mathsf{p}_{1l},\mathsf{q}_{1l})\in\mathsf{S}_1^{[2]}$,
$Z_{1l}^{[1]}$ is initialized to a corresponding part of $I_n$, while
$\mathsf{F}_{1l}^{[1]}$, $\mathsf{G}_{1l}^{[1]}$, and
$\mathsf{Z}_{1l}^{[1]}$ are undefined.
\subsubsection{Processing the block pivots}\label{sss:6.4.2}
In a step $\mathsf{j}$ (of a block sweep $\mathsf{s}$, which index we
omit from the superscripts of the matrices for simplicity when it is
implied by the context), the $l$th thread shortens its $J$ to
$\widehat{J}_{\mathsf{j}l}$ and its block column pairs from
$F_{\mathsf{j}l}^{}$ to $\widehat{F}_{\mathsf{j}l}^{}$ such that
$\widehat{H}_{\mathsf{j}l} \assgn \widehat{F}_{\mathsf{j}l}^{\ast} \widehat{J}_{\mathsf{j}l}^{} \widehat{F}_{\mathsf{j}l}^{} = F_{\mathsf{j}l}^{\ast} J^{} F_{\mathsf{j}l}^{}$,
and $G_{\mathsf{j}l}^{}$ to $\widehat{G}_{\mathsf{j}l}^{}$ such that
$\widehat{S}_{\mathsf{j}l} \assgn \widehat{G}_{\mathsf{j}l}^{\ast} \widehat{G}_{\mathsf{j}l}^{} = G_{\mathsf{j}l}^{\ast} G_{\mathsf{j}l}^{}$,
where $\widehat{F}_{\mathsf{j}l}^{}$ and
$\widehat{G}_{\mathsf{j}l}^{}$ are both square, of order $2\mathsf{w}$
or $2(\mathsf{w}-1)$ (the order in between the two is necessarily odd,
so the bordering from subsection~\ref{sss:6.3.2} should then be
applied).

There are two ways to shorten the block column pairs.  The more
efficient one is to form $\widehat{H}_{\mathsf{j}l}$ and
$\widehat{S}_{\mathsf{j}l}$ explicitly.  To do that,
$F_{\mathsf{j}l}^{}$ is copied to $\mathsf{F}_{\mathsf{j}l}^{}$, and
then the rows of $\mathsf{F}_{\mathsf{j}l}^{}$ are scaled by $J$.
Here, it is beneficial to have $J$ in the compact, run-length encoded
form.  Then, a single (parallel, or in our case, sequential)
\texttt{ZGEMM} call computes
$\widehat{H}_{\mathsf{j}l} = F_{\mathsf{j}l}^{\ast} (J^{} \mathsf{F}_{\mathsf{j}l}^{})$
and stores it temporarily into $\widehat{G}_{\mathsf{j}l}^{}$.  Then,
$\widehat{F}_{\mathsf{j}l}^{\ast}$ and $\widehat{J}_{\mathsf{j}l}$
(partitioned to the positive and the negative sign blocks) are
obtained by the Hermitian indefinite factorization with complete
pivoting (see section~\ref{s:4}),
$\widehat{H}_{\mathsf{j}l} = \widehat{F}_{\mathsf{j}l}^{\ast} \widehat{J}_{\mathsf{j}l}^{} \widehat{F}_{\mathsf{j}l}^{}$.
The factorization should reveal if $\widehat{H}_{\mathsf{j}l}$ is rank
deficient, in which case the process stops (the computation could be
retried with the (J)QR approach, as below, if the input data has been
preserved).  Finally, $\widehat{F}_{\mathsf{j}l}^{\ast}$ is copied,
with the transposition and the complex conjugation applied, to
$\widehat{F}_{\mathsf{j}l}^{}$.

A simpler procedure is used for $\widehat{G}_{\mathsf{j}l}$.  It
suffices to compute $G_{\mathsf{j}l}^{\ast} G_{\mathsf{j}l}^{}$ by a
single (parallel, or in our case, sequential) \texttt{ZHERK} call,
store $\widehat{S}_{\mathsf{j}l}$ temporarily to
$\widehat{Z}_{\mathsf{j}l}$, and perform either the diagonally-pivoted
Cholesky factorization, or even the Hermitian indefinite factorization
with complete pivoting, to obtain
$\widehat{S}_{\mathsf{j}l} = \widehat{G}_{\mathsf{j}l}^{\ast} \widehat{G}_{\mathsf{j}l}^{}$.
If $\widehat{S}_{\mathsf{j}l}$ is rank deficient or indefinite, the
algorithm stops with an error message.  A similar treatment has been
implemented for $\widehat{H}_{\mathsf{j}l}$, if $J = I$.  Otherwise,
$\widehat{G}_{\mathsf{j}l}^{\ast}$ is copied, with the transposition
and the complex conjugation applied, to
$\widehat{G}_{\mathsf{j}l}^{}$.  A \texttt{ZLASET} call finally
initializes $\widehat{Z}_{\mathsf{j}l}$ to $I$.

The (J)QR approach, similar to Phase~2, should be more accurate and
even necessary when the conditions of the matrices
$\widehat{H}_{\mathsf{j}l}$ and/or $\widehat{S}_{\mathsf{j}l}$ are so
large that factorizing them after forming them explicitly might
fail.  Both $F_{\mathsf{j}l}$ and $G_{\mathsf{j}l}$ should then be
copied to $\mathsf{F}_{\mathsf{j}l}$ and $\mathsf{G}_{\mathsf{j}l}$,
respectively, and the JQR factorization on $\mathsf{F}_{\mathsf{j}l}$,
followed by the column prepermutation and the column-pivoted QR
factorization on $\mathsf{G}_{\mathsf{j}l}$, should be performed to
obtain $\widehat{F}_{\mathsf{j}l}$ with $\widehat{J}_{\mathsf{j}l}$,
and $\widehat{G}_{\mathsf{j}l}$, respectively.  This has not been
implemented, though.

The block pivots $\widehat{F}_{\mathsf{j}l}$ with
$\widehat{J}_{\mathsf{j}l}$ and $\widehat{G}_{\mathsf{j}l}$ are handed
over to a version of the Level~1 (VP) algorithm to be transformed.
This, single-threaded VP version executes in the contexts of the
already running threads.  The Level~1 algorithm can either fully
(implicitly) diagonalize $\widehat{H}_{\mathsf{j}l}$ and
$\widehat{S}_{\mathsf{j}l}$, or at least iterate until a reasonably
high number of the inner sweeps has been attained
($\mathsf{C}_{\max} = 30$), in which case we talk about the Full Block
(FB) variant; or pass over the block pivots only a prescribed number
of times, e.g., once ($\mathsf{C}_{\max} = 1$), in the Block-Oriented
(BO) variant (for more details on both,
see~\cite{Novakovic-SingerSanja-SingerSasa-2015}).  The former variant
corresponds to a full two-sided annihilation of the off-diagonal of
$\widehat{H}_{\mathsf{j}l}$ and $\widehat{S}_{\mathsf{j}l}$, while the
latter implicitly reduces their off-diagonal norms.

Also, the BO variant exhibits similar execution times for every call of
the Level~1 routine across all threads in a block step, while those
can vary significantly, due to data, among the threads in the FB
variant.  On some platforms (e.g., the GPUs), that does not pose a
huge problem~\cite{Novakovic-2015}, but on the CPUs it can cause
delays on the synchronization primitives (OpenMP barriers) required
between the block steps~\cite{Novakovic-SingerSanja-SingerSasa-2015}.

In any case, the transformations applied in the Level~1 are
accumulated in $\widehat{Z}_{\mathsf{j}l}$, while the counters of all
and of ``big'' transformations are added atomically to the respective
counters (shared among the threads) for the current block sweep.  With
distributed memory, such counters can be updated by
\texttt{MPI\_Allreduce} collective calls.
\subsubsection{Updating and exchanging the block columns}\label{sss:6.4.3}
If no transformations have been applied in the Level~1,
$F_{\mathsf{j}l}$, $G_{\mathsf{j}l}$, and $Z_{\mathsf{j}l}$ are then
copied to $\mathsf{F}_{\mathsf{j}l}$, $\mathsf{G}_{\mathsf{j}l}$, and
$\mathsf{Z}_{\mathsf{j}l}$, respectively.  Otherwise,
$\mathsf{F}_{\mathsf{j}l} = F_{\mathsf{j}l} \widehat{Z}_{\mathsf{j}l}$,
$\mathsf{G}_{\mathsf{j}l} = G_{\mathsf{j}l} \widehat{Z}_{\mathsf{j}l}$, and
$\mathsf{Z}_{\mathsf{j}l} = Z_{\mathsf{j}l} \widehat{Z}_{\mathsf{j}l}$
(three \texttt{ZGEMM} calls).

All computation have thus far been local to a thread, apart from the
atomic operations.  Now, by looking to
$\mathsf{S}_{\mathsf{j}+1}^{[2]}$, it is easy to figure out which of
the two block columns should be retained by the thread, and to which
thread the other block column has to be sent (in fact, the whole
communication pattern has been precomputed).  Also, the block columns
can swap roles: the first one in the current step can become the
second one, in the presently owning or in the receiving thread, and
vice versa, in the next step.  A thread might even send away both its
block columns (each to a different recipient), and receive two new
block columns (each from a different sender).

It might seem unnecessary to perform the physical exchanges of data on
the shared memory, but the reason behind them is twofold.  First, most
modern machines have their memory paritioned according to the speed of
access by, or ``proximity'' to, each CPU (NUMA), and it is beneficial
to bring the data close to (a thread bound to) the CPU that processes
it.  Second, such a design makes the algorithm convertible to a
distributed memory one (e.g., by using an MPI process for what a
thread does now).

To perform the block column copies, a thread has to know the memory
addresses of the storage of the threads it communicates to.  All such
addresses are kept available to all threads, and the block columns
from the ``shadow'' storage are copied by their present owners to the
regular storage of their future ones.  Before a block column pair can
be copied (from or to), it has to be updated first, so there is an
OpenMP barrier between the three updates above and the three copying
actions (by two \texttt{ZLACPY} calls each, one per a block column).
Also, a thread cannot continue with the next block step until it has
the new data ready and its shadow storage available, which is
enforced for all threads simultaneously by placing another barrier
after the copying actions.

It is a legitimate question if it is worthwhile (and in what
circumstances) to try hiding the communication behind the computation,
maybe by relying on some tasking mechanism.  For example,
$\mathsf{F}_{\mathsf{p}_{\mathsf{j}l}}$ can be copied to
$F_{\mathsf{q}_{\mathsf{j}+1,l'}}$, and
$\mathsf{F}_{\mathsf{q}_{\mathsf{j}l}}$ to
$F_{\mathsf{p}_{\mathsf{j}+1,l''}}$, while $G_{\mathsf{j}l}$,
$G_{\mathsf{j}l'}$, and $G_{\mathsf{j}l''}$ are being updated.  We
leave those considerations for a future work.

\paragraph{Completing a sweep}%
At the end of a block sweep, the transformation counters are read and
reset.  If no ``big'' transformations have been applied in any of the
threads in the block sweep, the process stops, and the outputs, now
including $\Lambda_F$, $\Lambda_G$, and $\Lambda$, are generated,
piecewise per thread, as described in subsection~\ref{sss:6.1.7} and
section~\ref{s:3}.
\subsection{Testing}\label{ss:6.5}
The aim of testing was to establish what algorithm variant to
recommend for practice, as well as should Phase~2 be employed, and if
so, when.
\subsubsection{Blocked vs.\ pointwise algorithms}\label{sss:6.5.1}
In Table~\ref{tbl:L1BOFB} it is shown that the Level~2 (BO) algorithm
is several times faster than both the Level~1 (VP) and the Level~2
(FB) algorithms, and that blocking in general gives a significant
advantage over the pointwise approach.  Also, having more threads, and
thus the smaller block pivot orders, benefits the FB algorithm, since
the speedup with 64 versus 32 threads is there more than twofold.  On
the contrary, the speedup with twice more threads is less than twofold
for the BO algorithm, since the formation of the block pivots in the
BO variant takes a bigger portion of the overall time, compared to the
Level~1 inner iterations.  Nevertheless, the Level~2 (BO) algorithm is
our choice for Phase~3.
\begin{table}[hbt]
  {\footnotesize\caption{The wall execution time (wtime), with 32, 64,
      and \emph{48} threads, of the three algorithm candidates.}
    \label{tbl:L1BOFB}
    \begin{center}
      \begin{tabular}{@{}cccccccc@{}}
        \toprule
        \multirow[c]{2}{*}{ID} & \multicolumn{2}{c}{Level~1 (VP) wtime [s]} & \multicolumn{2}{c}{Level~2 (FB) wtime [s]} & \multicolumn{3}{c}{Level~2 (BO) wtime [s]}\\
        & 32 thr. & 64 thr. & 32 thr. & 64 thr. & 32 thr. & 64 thr. & \emph{48 thr\/}.\\
        \toprule
        A1 & $\hphantom{00}828.58$ & $\hphantom{00}836.50$ & $\hphantom{00}372.00$ & $\hphantom{00}222.81$ & $\hphantom{0}125.38$ & $\hphantom{0}128.51$ & $\hphantom{00}39.39$\\
        \midrule
        A2 & $\hphantom{0}4528.13$ & $\hphantom{0}4533.41$ & $\hphantom{0}2215.07$ & $\hphantom{00}914.32$ & $\hphantom{0}634.72$ & $\hphantom{0}481.92$ & $\hphantom{0}169.97$\\
        \midrule
        A3 &            $19094.72$ &            $19357.67$ & $\hphantom{0}9253.03$ & $\hphantom{0}3539.32$ &            $2392.01$ &            $1474.61$ & $\hphantom{0}569.85$\\
        \midrule
        A4 &            $58910.73$ &            $54938.04$ &            $30925.99$ &            $11449.48$ &            $8003.41$ &            $4528.69$ & $2128.77$\\
        \midrule
        B1 & $\hphantom{00}254.69$ & $\hphantom{00}266.73$ & $\hphantom{00}135.87$ & $\hphantom{00}101.21$ & $\hphantom{00}46.36$ & $\hphantom{00}50.43$ & $\hphantom{00}16.91$\\
        \midrule
        B2 & $\hphantom{0}1379.01$ & $\hphantom{0}1461.47$ & $\hphantom{00}619.19$ & $\hphantom{00}328.24$ & $\hphantom{0}188.76$ & $\hphantom{0}173.86$ & $\hphantom{00}53.28$\\
        \midrule
        B3 & $\hphantom{0}6377.36$ & $\hphantom{0}6249.12$ & $\hphantom{0}2870.85$ & $\hphantom{0}1142.89$ & $\hphantom{0}789.13$ & $\hphantom{0}529.44$ & $\hphantom{0}190.45$\\
        \midrule
        B4 &            $22002.45$ &            $22421.95$ & $\hphantom{0}9620.59$ & $\hphantom{0}4015.04$ &            $2361.69$ &            $1506.33$ & $\hphantom{0}579.38$\\
        \bottomrule
      \end{tabular}
  \end{center}}
\end{table}
\subsubsection{Phase2 benefits}\label{sss:6.5.2}
Phase~3 might have been run directly on the $\widetilde{F}$,
$\widetilde{J}$, and $\widetilde{G}$.  In Table~\ref{tbl:SQvsTS} it is
shown that the preprocessing of the tall-and-skinny inputs by Phase~2
into the square ones for Phase~3 is preferred, both time-wise and
sweep-wise, but that advantage diminishes as the inputs approach a
square-like form ($m\gtrsim n$).
\begin{table}[hbt]
  {\footnotesize\caption{The wall execution time (wtime) with 32, 64,
      and \emph{48} threads, the speedups, and the numbers of block
      sweeps of Phase~3 on the tall-and-skinny and the square inputs
      (obtained by Phase~2).}
    \label{tbl:SQvsTS}
    \begin{center}
      \begin{tabular}{@{}cccccccc@{}}
        \toprule
        \multirow[c]{2}{*}{ID} & \multicolumn{2}{c}{BO tall-and-skinny $(\diamond)$} & \multicolumn{3}{c}{Phase~2 \&\ BO square $(\star)$} & \multicolumn{2}{c}{speedup $(\diamond)/(\star)$ \&\ sweeps}\\
        & 32 thr. & 64 thr. & 32 thr. & 64 thr. & \emph{48 thr\/}. & 32 thr. & 64 thr.\\
        \toprule
        A1 & $\hphantom{00}590.92$ &            $1035.61$ & $\hphantom{00}298.71$ & $\hphantom{0}315.04$ & $\hphantom{0}107.73$ & $1.978; 15\vert 14$ & $3.287; 15\vert 14$\\
        \midrule
        A2 & $\hphantom{0}1813.66$ &            $1600.17$ & $\hphantom{0}1098.79$ & $\hphantom{0}973.41$ & $\hphantom{0}358.56$ & $1.651; 17\vert 16$ & $1.644; 16\vert 16$\\
        \midrule
        A3 & $\hphantom{0}5393.96$ &            $3457.28$ & $\hphantom{0}3457.09$ &            $2594.89$ & $\hphantom{0}961.61$ & $1.560; 19\vert 17$ & $1.332; 17\vert 16$\\
        \midrule
        A4 &            $13293.75$ &            $7734.08$ &            $10185.98$ &            $6808.09$ &            $2899.21$ & $1.305; 19\vert 18$ & $1.136; 18\vert 18$\\
        \midrule
        B1 & $\hphantom{00}609.11$ & $\hphantom{0}992.56$ & $\hphantom{00}239.55$ & $\hphantom{0}240.14$ & $\hphantom{00}85.10$ & $2.543; 12\vert 13$ & $4.133; 12\vert 12$\\
        \midrule
        B2 & $\hphantom{0}1665.46$ &            $1716.25$ & $\hphantom{00}725.15$ & $\hphantom{0}680.55$ & $\hphantom{0}244.13$ & $2.297; 14\vert 13$ & $2.522; 13\vert 13$\\
        \midrule
        B3 & $\hphantom{0}3595.73$ &            $3126.35$ & $\hphantom{0}2091.06$ &            $1733.91$ & $\hphantom{0}620.94$ & $1.720; 15\vert 15$ & $1.803; 14\vert 14$\\
        \midrule
        B4 & $\hphantom{0}8031.08$ &            $6147.36$ & $\hphantom{0}5194.34$ &            $3981.79$ &            $1440.84$ & $1.546; 15\vert 15$ & $1.544; 15\vert 15$\\
        \bottomrule
      \end{tabular}
  \end{center}}
\end{table}
\subsubsection{Generalized eigenvalues}\label{sss:6.5.3}
Please see the Figures~S.2 and S.3 in the supplementary material,
depicting the generalized eigenvalues
$\Lambda(\text{A1})$--$\Lambda(\text{A4})$ and
$\Lambda(\text{B1})$--$\Lambda(\text{B4})$, computed with 64 threads
by the Level~2 (BO) Phase~3 GHSVD after shortening in Phase~2.
The eigenvalues obtained with 32 threads differ a few ulps at most.
\subsubsection{Comparison with \texttt{ZHEGV}(\texttt{D})}\label{sss:6.5.4}
In Table~\ref{tbl:ZHEGVD64} the wall times for the explicit formation
of $H$ (by the $\widetilde{J}$-scaling of $\widetilde{F}$ and the
\texttt{ZGEMM} matrix multiplication) and $S$ (by the \texttt{ZHERK}
matrix multiplication), as well as for the LAPACK's generalized
Hermitian eigensolvers \texttt{ZHEGV} and \texttt{ZHEGVD} on $(H,S)$,
left in MCDRAM when possible, are shown, alongside the speedups of
this approach (i.e., of forming of $H$ and $S$ and then calling either
\texttt{ZHEGV} or \texttt{ZHEGVD}, with the eigenvectors also
computed) vs.~the Level~2 (BO) Phase~3 GHSVD on the inputs shortened
by Phase~2, with 64 and \emph{48\/} threads.
\begin{table}[hbt]
  {\footnotesize\caption{The wall times for the explicit formation of
      $H$ and $S$, combined with those for a LAPACK's generalized
      Hermitian eigensolver (\texttt{ZHEGV} or  \texttt{ZHEGVD}), and
      the speedup vs.~$(\star)$, with 64 and \emph{48} threads.}
    \label{tbl:ZHEGVD64}
    \begin{center}
      \setlength{\tabcolsep}{4.5pt}
      \begin{tabular}{@{}cccccccccc@{}}
        \toprule
        \multirow[c]{2}{*}{ID} & \#\ of & \multicolumn{4}{c}{(max.~of 2 runs for $H,S$) wall time [s] for} & \multicolumn{2}{c}{total wall time [s] with} & \multicolumn{2}{c}{speedup vs.~$(\star)$}\\
        & thrs. & $H=\widetilde{F}^{\ast}\widetilde{J}\widetilde{F}$ & $S=\widetilde{G}^{\ast}\widetilde{G}$ & \texttt{ZHEGV} & \texttt{ZHEGVD} & \texttt{ZHEGV}~$(\triangleleft)$ & \texttt{ZHEGVD}~$(\triangleright)$ & $(\star)/(\triangleleft)$ & $(\star)/(\triangleright)$\\
        \toprule
        \multirow[c]{2}{*}{A1} & 64 & $\hphantom{0}1.490$ & $\hphantom{0}1.519$ & $\hphantom{0}11.896$ & $\hphantom{00}5.573$ & $\hphantom{0}14.905$ & $\hphantom{00}8.317$ & $21.137$ & $37.881$\\
                        & \emph{48} & $\hphantom{0}1.316$ & $\hphantom{0}0.536$ & $\hphantom{00}2.875$ & $\hphantom{00}1.800$ & $\hphantom{00}4.727$ & $\hphantom{00}3.185$ & $22.793$ & $33.825$\\
        \midrule
        \multirow[c]{2}{*}{A2} & 64 & $\hphantom{0}3.959$ & $\hphantom{0}3.323$ & $\hphantom{0}38.197$ & $\hphantom{0}18.890$ & $\hphantom{0}45.479$ & $\hphantom{0}25.913$ & $21.403$ & $37.565$\\
                        & \emph{48} & $\hphantom{0}2.368$ & $\hphantom{0}1.277$ & $\hphantom{0}12.722$ & $\hphantom{00}8.504$ & $\hphantom{0}16.349$ & $\hphantom{0}12.148$ & $21.931$ & $29.515$\\
        \midrule
        \multirow[c]{2}{*}{A3} & 64 &            $11.642$ & $\hphantom{0}7.696$ &            $100.320$ & $\hphantom{0}56.839$ &            $119.448$ & $\hphantom{0}75.966$ & $21.724$ & $34.158$\\
                        & \emph{48} & $\hphantom{0}5.808$ & $\hphantom{0}3.201$ & $\hphantom{0}53.464$ & $\hphantom{0}37.089$ & $\hphantom{0}62.472$ & $\hphantom{0}46.098$ & $15.393$ & $20.860$\\
        \midrule
        \multirow[c]{2}{*}{A4} & 64 &            $24.399$ &            $15.369$ &            $284.750$ &            $168.956$ &            $324.281$ &            $208.513$ & $20.994$ & $32.651$\\
                        & \emph{48} &            $12.692$ & $\hphantom{0}7.095$ &            $160.374$ &            $119.576$ &            $180.157$ &            $139.290$ & $16.093$ & $20.814$\\
        \midrule
        \multirow[c]{2}{*}{B1} & 64 & $\hphantom{0}1.348$ & $\hphantom{0}1.909$ & $\hphantom{00}5.399$ & $\hphantom{00}2.573$ & $\hphantom{00}8.656$ & $\hphantom{00}5.575$ & $27.742$ & $43.076$\\
                        & \emph{48} & $\hphantom{0}0.866$ & $\hphantom{0}0.507$ & $\hphantom{00}1.055$ & $\hphantom{00}0.603$ & $\hphantom{00}2.427$ & $\hphantom{00}1.954$ & $35.058$ & $43.550$\\
        \midrule
        \multirow[c]{2}{*}{B2} & 64 & $\hphantom{0}4.401$ & $\hphantom{0}3.786$ & $\hphantom{0}17.015$ & $\hphantom{00}8.750$ & $\hphantom{0}24.987$ & $\hphantom{0}16.648$ & $27.236$ & $40.879$\\
                        & \emph{48} & $\hphantom{0}2.270$ & $\hphantom{0}1.281$ & $\hphantom{00}4.284$ & $\hphantom{00}2.791$ & $\hphantom{00}7.799$ & $\hphantom{00}6.342$ & $31.302$ & $38.493$\\
        \midrule
        \multirow[c]{2}{*}{B3} & 64 &            $11.213$ &            $10.047$ & $\hphantom{0}45.139$ & $\hphantom{0}24.654$ & $\hphantom{0}63.692$ & $\hphantom{0}45.914$ & $27.223$ & $37.764$\\
                        & \emph{48} & $\hphantom{0}5.508$ & $\hphantom{0}3.104$ & $\hphantom{0}15.564$ & $\hphantom{0}12.399$ & $\hphantom{0}24.162$ & $\hphantom{0}20.904$ & $25.699$ & $29.704$\\
        \midrule
        \multirow[c]{2}{*}{B4} & 64 &            $23.092$ &            $16.446$ &            $114.787$ & $\hphantom{0}62.469$ &            $150.798$ &            $102.006$ & $26.405$ & $39.035$\\
                        & \emph{48} &            $12.003$ & $\hphantom{0}6.790$ & $\hphantom{0}51.817$ & $\hphantom{0}41.265$ & $\hphantom{0}70.570$ & $\hphantom{0}59.836$ & $20.417$ & $24.080$\\
        \bottomrule
      \end{tabular}
  \end{center}}
\end{table}
Please see Table~S.2 in the supplementary material for the results
with 32 and \emph{24\/} threads.

The results suggest that both the CPU's clock and the size of its
low-level caches play a crucial role in performance of our approach
for the problems large enough.

\looseness=-1
The main advantage of our approach is its ability to compute the
generalized eigenproblem accurately when the LAPACK-based one can
fail, not least due to the squaring of the condition number
$\mathop{\kappa_2}(\widetilde{G})$ in the forming of
$S=\widetilde{G}^{\ast}\widetilde{G}$, i.e.,
\begin{displaymath}
  \kappa_2^{}(S)=\lambda^{}_{\max}(S)/\lambda^{}_{\min}(S)=\sigma^2_{\max}(\widetilde{G})/\sigma^2_{\min}(\widetilde{G})=\kappa_2^2(\widetilde{G}),
\end{displaymath}
that consequently leads to a failure in the Cholesky factorization, or
to an unacceptable inaccuracy in the generalized eigenvalues.  We
believe that the demonstrated slowdown is a reasonable tradeoff for a
reliable backup alternative in those difficult cases.
%
%
\section{Phase~4 -- optional computation of the full G(H)SVD}\label{s:7}
%
%
In Phase~4 the right generalized singular vector matrix
$X^{} = Z^{-1}$ is computed by the LU factorization with complete
pivoting, $Z^{} = P^T L^{} U^{} Q^T$, from the LAPACK routine
\texttt{ZGETC2}, followed by solving a linear system $Z X = I_n$ for
$X$, with $Z$ factored as above, by calling the sequential routine
\texttt{ZGESC2} in a parallel do loop.  Each of $n$ loop iterations
solves for one column of $X$, and the iteration space is divided among
the same number of OpenMP threads, also available to \texttt{ZGETC2},
as used for the previous phases.

Should the G(H)SVD of the tall-and-skinny factors be required,
$\widetilde{Z}^{} = P_2^T Z^{}$ and
$\widetilde{X}^{} = \widetilde{Z}^{-1}$ could also be computed in
Phase~4, together with $\widetilde{U} = \widetilde{Q}_F^{} U$ and
$\widetilde{V} = \widetilde{Q}_G^{} V$.
\subsection{Relative errors in the full GHSVD}\label{ss:7.1}
To measure the accuracy of the Phase~3 algorithm, we can look at the
Frobenius norm of the error in the obtained GHSVD decomposition,
relative to the Frobenius norm of the original matrix, i.e.,
\begin{equation}
  \|F - U \Sigma_F X\|_F / \|F\|_F
  \quad\text{and}\quad
  \|G - V \Sigma_G X\|_F / \|G\|_F.
  \label{7.1}
\end{equation}

For a detailed description of computing the relative errors, see the
supplementary material, subsection~S.3.1, and for the full accuracy
results, see Table~S.3 therein.

The relative errors in the decomposition of $F$ range from
$1.25\cdot 10^{-13}$ to $3.52\cdot 10^{-13}$, and from
$1.12\cdot 10^{-13}$ to $7.22\cdot 10^{-13}$, for the datasets A and
B, respectively.  In the decomposition of $G$, the relative errors are
slightly lower, ranging from $9.47\cdot 10^{-14}$ to
$8.23\cdot 10^{-13}$, and from $7.98\cdot 10^{-14}$ to
$4.73\cdot 10^{-13}$ for the datasets A and B, respectively.  The
relative errors are similar regardless the number of threads used in
Phase~3.

From the range of errors it can be concluded both that our datasets
are not highly ill-conditioned, and that in those cases the Phases~3
and 4 do not behave erratically.  A rigorous stability analysis of the
GHSVD method remains open for the future work.
\subsection{Comparison with \texttt{ZGGSVD3}}\label{ss:7.2}
The Phase~3 algorithm can also be used for the ordinary GSVD by
setting $J = I$.  It is therefore reasonable to compare it to the
\texttt{ZGGSVD3} LAPACK routine (from the parallel Intel MKL) serving
the same purpose.
\subsubsection{Dataset}\label{sss:7.2.1}
A dataset C, comprising five Hermitian matrix pairs, has been
generated by a call to the \texttt{ZLATMS} LAPACK testing routine for
each matrix of the full bandwidth, with its pseudorandom eigenvalues
uniformly distributed in $(0,1)$.  The matrix orders are
$k\cdot 1000$, $1\le k\le 5$, and the full GSVD is required.
\subsubsection{Timing results}\label{sss:7.2.2}
In Table~\ref{tbl:gsvd64} the speedup of Phase~3 followed by Phase~4,
both with 64 threads, vs.\ \texttt{ZGGSVD3} is shown.  The speedup
with 32 threads reaches a lower peak for $n=5000$, 140 times, as shown
in the supplementary material, Table~S.4.
\begin{table}[hbt]
  {\footnotesize\caption{The wall execution time (wtime) and the
      speedup of the Phases~3 (with $J = I$) and 4 versus the
      \texttt{ZGGSVD3} LAPACK routine on the set C, with 64 threads
      and $n$ denoting the order of the matrices.}
    \label{tbl:gsvd64}
    \begin{center}
      \begin{tabular}{@{}cccccc@{}}
        \toprule
        \multirow[c]{2}{*}{$n$} & \texttt{ZGGSVD3} & Phase~3 & Phase~4 & Phases~3~\&~4 & speedup \\
        & wtime [s] $(\bullet)$ & wtime [s] \&\ sweeps & wtime [s] & wtime [s] $(\circ)$ & $(\bullet)/(\circ)$ \\
        \toprule
        $1000$ & $\hphantom{000}399.47$ & $\hphantom{00}8.32; 13$ & $\hphantom{00}2.84$ & $\hphantom{0}11.16$ & $\hphantom{0}35.78$ \\
        \midrule
        $2000$ & $\hphantom{00}5935.03$ & $\hphantom{0}35.52; 14$ & $\hphantom{0}20.28$ & $\hphantom{0}55.80$ &            $106.37$ \\
        \midrule
        $3000$ & $\hphantom{0}21880.11$ &            $100.71; 16$ & $\hphantom{0}69.48$ &            $170.19$ &            $128.57$ \\
        \midrule
        $4000$ & $\hphantom{0}54233.01$ &            $191.32; 16$ &            $180.70$ &            $372.02$ &            $145.79$ \\
        \midrule
        $5000$ &            $107424.26$ &            $332.29; 17$ &            $327.97$ &            $660.26$ &            $162.70$ \\
        \bottomrule
      \end{tabular}
  \end{center}}
\end{table}
From the tests it is evident that Phase~3 (with Phase~4 if needed)
forms a very competitive, highly parallel algorithm for the ordinary
GSVD as well, while being adaptable to a wide variety of the modern
high performance computing hardware.
%
%
\section{Future work}\label{s:8}
%
%
A (multi-)GPU version of the implicit Hari--Zimmermann algorithm for
the ordinary GSVD~\cite{Novakovic-2017,Novakovic-Singer-2019} shows
the promising results and complements the CPU implementation described
herein.  Therefore, a GPU implementation of the GHSVD and Phase~1
might be a practical companion to the present research.
%
%
\section*{Acknowledgments}
%
%
The authors would like thank the anonymous referees for their
suggestions on improving the contents and the presentation of this
manuscript.
%
%

%
%
\end{document}